\documentclass[reqno]{amsart}

\usepackage{amssymb}
\usepackage{amsmath}
\usepackage{bbm}
\usepackage{mathrsfs}
\usepackage{xypic}
\usepackage{textcomp}

\makeatletter
\def\hsmash{\relax 
  \ifmmode\def\next{\mathpalette\mathhsm@sh}\else\let\next\makehsm@sh
  \fi\next}
\def\makehsm@sh#1{\setbox\z@\hbox{#1}\finhsm@sh}
\def\mathhsm@sh#1#2{\setbox\z@\hbox{$\m@th#1{#2}$}\finhsm@sh}
\def\finhsm@sh{\wd\z@\z@ \box\z@}
\makeatother

\newcommand{\br}{ }
\newcommand{\brr}{, }

\newtheorem{theorem}{Theorem}[section]
\newtheorem{theo}[theorem]{Theorem}
\newtheorem{mth}[theorem]{Theorem}
\newtheorem{lem}[theorem]{Lemma}
\newtheorem{prop}[theorem]{Proposition}

\theoremstyle{definition}
\newtheorem{ttt}[theorem]{}
\newtheorem{defi}[theorem]{Definition}
\newtheorem{defs}[theorem]{Definitions}
\newtheorem{nota}[theorem]{Notation}
\newtheorem{conv}[theorem]{Convention}
\newtheorem{cau}[theorem]{Caution}
\newtheorem{ex}[theorem]{Example}

\newtheorem{defio}[theorem]{Definition}
\newtheorem{exo}[theorem]{Example}
\newtheorem{algoo}[theorem]{Algorithm}

\theoremstyle{remark}

\newtheorem{rem}[theorem]{Remark}
\newtheorem{rems}[theorem]{Remarks}

\newtheorem{remo}[theorem]{Remark}
\newtheorem{remso}[theorem]{Remarks}

\numberwithin{equation}{section}

\newcommand{\calC}{\mathscr{C}}
\newcommand{\calE}{\mathscr{E}}
\newcommand{\calI}{\mathscr{I}}
\newcommand{\calL}{\mathscr{L}}
\newcommand{\calM}{\mathscr{M}}
\newcommand{\calO}{\mathscr{O}}
\newcommand{\calP}{\mathscr{P}}

\newcommand{\calU}{\mathscr{U}}
\newcommand{\calV}{\mathscr{V}}

\newcommand{\wcalM}{\,\widetilde{\!\calM}}

\newcommand{\bbQ}{\mathbbm{Q}}
\newcommand{\bbZ}{\mathbbm{Z}}

\newcommand{\frakp}{\mathfrak{p}}

\newcommand{\SL}{\mathop{\rm SL}\nolimits}
\newcommand{\GL}{\mathop{\rm GL}\nolimits}
\newcommand{\PGL}{\mathop{\rm PGL}\nolimits}

\newcommand{\im}{\mathop{\rm im}\nolimits}
\newcommand{\Aut}{\mathop{\rm Aut}\nolimits}
\newcommand{\Mor}{\mathop{\rm Mor}\nolimits}
\newcommand{\Spec}{\mathop{\rm Spec}\nolimits}

\newcommand{\tr}{\mathop{\rm tr}\nolimits}
\newcommand{\Tr}{\mathop{\rm Tr}\nolimits}
\newcommand{\N}{\mathop{\rm N}\nolimits}
\newcommand{\Sym}{\mathop{\rm Sym}\nolimits}

\newcommand{\res}{\mathop{\rm res}\nolimits}

\newcommand{\Hom}{{\rm Hom}}
\newcommand{\Pic}{{\rm Pic}}
\newcommand{\Gal}{{\rm Gal}}
\renewcommand{\div}{{\rm div}}
\newcommand{\pr}{{\rm pr}}
\newcommand{\Q}{{\rm Q}}
\newcommand{\id}{{\rm id}}
\newcommand{\sing}{{\rm sing}}
\newcommand{\Cl}{{\rm Cl}}

\newcommand{\bA}{{\bf A}}
\newcommand{\bP}{{\bf P}}

\newcommand{\pmodulo}[1]{\nobreak\ifinner\mkern8mu\else\mkern18mu\fi
 (\textup{mod}\,\,#1)}

\newcommand{\ratarrow}{$%
$\definemorphism{rat}\dashed\tip\notip%
\spreaddiagramcolumns{-12pt}%
 - \!\!\diagram%
\rrat &
\enddiagram\!\!$%
$}

\newcounter{abc}
\newenvironment{abc}{\begin{list}{\rm \alph{abc}) }%
{\usecounter{abc} \leftmargin=0.0pt \labelsep=0.0pt %
\listparindent=0.0pt \labelwidth=0.0pt \parsep=\smallskipamount %
\itemsep=0.0pt \topsep=0.0pt \partopsep=\smallskipamount}}{\end{list}}

\newcounter{iii}
\newenvironment{iii}{\begin{list}{\rm \roman{iii}) }%
{\usecounter{iii} \leftmargin=0.0pt \labelsep=0.0pt %
\listparindent=0.0pt \labelwidth=0.0pt \parsep=\smallskipamount%
 \itemsep=0.0pt \topsep=0.0pt \partopsep=\smallskipamount}}{\end{list}}

\def\rightend#1#2{{%
 \leavevmode\nobreak\hskip .5em plus 1fil
 \penalty600 \hskip 0pt plus -1filll
 \vadjust{}\nobreak\hskip 0pt plus 1filll%
 #1\parfillskip=#2\relax \par}}

\def\eop{\ifmmode\rule[-22pt]{0pt}{1pt}\ifinner\tag*{$\square$}\else\eqno{\square}\fi\else\rightend{$\square$}{0pt}\fi}

\begin{document}
\title[Inverse Galois problem for cubic surfaces]{Moduli spaces and the inverse Galois problem \\for cubic~surfaces}

\author{Andreas-Stephan Elsenhans}
\address{School of Mathematics and Statistics F07, University of Sydney, NSW 2006, Sydney, Australia, {\tt http://www.staff\!.\!uni-bayreuth.de/$\sim$bt270951/}}
\email{stephan@maths.usyd.edu.au}
\thanks{The first author was supported in part by the Deutsche Forschungsgemeinschaft (DFG) through a funded research~project.}

\author{J\"org Jahnel}
\address{D\'epartement Mathematik, Universit\"at Siegen, Walter-Flex-Str.~3, D-57068 Siegen, Germany, {\tt http://www.uni-math.gwdg.de/jahnel}}
\email{jahnel@mathematik.uni-siegen.de}
\thanks{}


\date{}
\dedicatory{}

\begin{abstract}
We~study the moduli space
$\smash{\wcalM}$
of marked cubic~surfaces. By~classical work of A.\,B.~Coble, this has a compactification
$\smash{\widetilde{M}}$,
which is linearly acted upon by the
group~$W(E_6)$.
$\smash{\widetilde{M}}$~is
given as the intersection of 30 cubics
in~$\bP^9$.
For~the morphism
$\smash{\wcalM \to \bP(1,2,3,4,5)}$
forgetting the marking, followed by Clebsch's invariant map, we give explicit~formulas. I.e.,~Clebsch's invariants are expressed in terms of Coble's irrational~invariants. As~an application, we give an affirmative answer to the inverse Galois problem for cubic surfaces
over~$\bbQ$.
\end{abstract}

\maketitle

\section*{Introduction}

Cubic~surfaces have been intensively studied by the geometers of the 19th~century. For~example, it was proven at that time that there are exactly 27~lines on every smooth cubic~surface. Further,~the configuration of the 27~lines is highly~symmetric. The~group of all permutations respecting the intersection pairing is isomorphic to the Weyl
group~$W(E_6)$
of
order~$51\,840$.

The~concept of a moduli scheme is by far more~recent. Nevertheless,~there are two kinds of moduli schemes for smooth cubic surfaces and both have their origins in classical invariant~theory.

On~one hand, there is the coarse moduli scheme
$\wcalM$
of smooth cubic~surfaces. This~scheme is essentially due to G.~Salmon~\cite{Sa} and A.~Clebsch~\cite{Cl}. In~fact, in a modern language, Clebsch's result from 1861 states that there is an open embedding
$\Cl\colon \calM \hookrightarrow \bP(1,2,3,4,5)$
into the weighted projective space of weights
$1, \ldots, 5$.

On~the other hand, one has the fine moduli scheme
$\smash{\wcalM}$
of smooth cubic surfaces with a marking on the 27~lines. The~marking plays the role of a rigidification and excludes all~automorphisms. That~is why a fine moduli scheme may~exist. It~has its origins in the work of A.~Cayley~\cite{Ca}. An~embedding
into~$\bP^9$
as an intersection of 30 cubics is due to A.\,B.~Coble~\cite{Co3} and dates back to the year~1917.

The~two moduli spaces are connected by the canonical, i.e.~forgetful, morphism
$\pr\colon \smash{\wcalM} \to \calM$.
This~is a finite flat morphism of
degree~$51\,840$.
Its~ramification locus corresponds exactly to the cubic surfaces having nontrivial~automorphisms.\bigskip

\noindent
{\bf Explicit formulas.}
In~Theorem~\ref{Formeln}, we will give an explicit description of
$\pr\colon \smash{\wcalM} \to \calM$.
In~other words, given a smooth cubic
surface~$C$
with a marking on its 27~lines, we provide explicit formulas expressing Clebsch's invariants
of~$C$
in terms of Coble's, so-called irrational,~invariants. From~a formal point of view, this result seems to be~new.

But~there can be no doubt that its essence, the existence of such formulas, has been clear to A.\,Coble, as~well. Only~due to the lack of computers, they could not be worked out at the time, with the exception of the very~first. In~fact, our approach is a combination of classical invariant theory with modern computer~algebra.\bigskip

\noindent
{\bf A solution to the equation problem.}
As~the main result of the article, we consider the following application of Theorem~\ref{Formeln}. Given~an abstract point on the moduli space of marked cubic surfaces, we deliver an algorithm that produces a concrete cubic surface from~it.

This~algorithm is a combination of the explicit formulas for
$\pr\colon \smash{\wcalM} \to \calM$
with an algorithmic solution to the so-called equation problem for cubic~surfaces. I.e.,~the problem to determine a concrete cubic surface from a given value of Clebsch's invariant~vector. This~was seemingly considered hopeless for a long time, but, today, it essentially comes down to the explicit computation of a Galois descent, cf.~\ref{sym} and Algorithm~\ref{expl}.\bigskip

\noindent
{\bf A further application.}
When~$C$
is a cubic surface
over~$\bbQ$,
the absolute Galois group
$\Gal(\overline\bbQ/\bbQ)$
operates on the 27~lines. This~means, after having fixed a marking on the lines, there is a
homomorphism~$\rho\colon \Gal(\overline\bbQ/\bbQ) \to W(E_6)$.
One~says that the Galois
group~$\Gal(\overline\bbQ/\bbQ)$
acts upon the lines
of~$C$
via~$G := \im \rho\subseteq W(E_6)$.
When~no marking is chosen, the
subgroup~$G$
is determined only up to~conjugation.

As~an application of the considerations on moduli schemes, we obtain the following affirmative answer to the inverse Galois problem for smooth cubic surfaces
over~$\bbQ$.

\begin{mth}\medskip
Let\/~$\mathfrak{g}$
be an arbitrary conjugacy class of subgroups
of\/~$W(E_6)$.
Then~there exists a smooth cubic
surface\/~$C$
over\/~$\bbQ$
such that the Galois group acts upon the lines
of\/~$C$
via a subgroup\/
$G \subseteq W(E_6)$
belonging to the conjugacy
class\/~$\mathfrak{g}$.
\end{mth}\medskip

\noindent
The~fundamental idea of the proof is as follows. We~describe a twist
$\smash{\wcalM_{\!\!\rho}}$
of
$\smash{\wcalM}$,
representing cubic surfaces with a marking that is acted upon by the absolute Galois group via a prescribed homomorphism
$\rho\colon \Gal(\overline\bbQ/\bbQ) \to W(E_6)$.
The~\mbox{$\bbQ$-rational}
points on this scheme correspond to the cubic surfaces of the type sought~for.

We~do not have the universal family at our disposal, a least not in a sufficiently explicit~form. Thus,~we calculate Clebsch's invariants of the cubic surface from the projective coordinates of the point found, i.e.~from the irrational invariants of the cubic~surface. Finally,~we recover the surface from Clebsch's invariants.\medskip

\noindent
{\bf The list.}
The~complete list of our examples is available at both author's web pages as a file named {\tt kub\_fl\_letzter\_teil.txt}. The~numbering of the conjugacy classes we use is that produced by {\tt gap}, version 4.4.12. This~numbering is reproducible, at least in our version of~{\tt gap}. It~coincides with the numbering used in our previous~articles.\medskip

\noindent
{\bf Acknowledgement.}
We~wish to thank Claus Fieker for several hints on how to handle in~{\tt magma} the relatively large number fields that occurred in the computations related to this~project.

\section{The moduli scheme of marked cubic surfaces}

The purpose of this section is mainly to fix notation and to recall some results that are more or less known.

\begin{defs}\smallskip
\begin{iii}
\item
Let~$S$
be any~scheme. Then,~by a {\em family of cubic surfaces\/}
over~$S$
or simply a {\em cubic surface\/}
over~$S$,
we mean a flat morphism
$p \colon C \to S$
such that there exist a
\mbox{rank-$4$}
vector bundle
$\calE$
on~$S$,
a non-zero section
$c \in \Gamma(\calO(3), \bP(\calE))$,
and an isomorphism
$\smash{\div(c) \stackrel{\cong}{\longrightarrow} C}$
of
$S$-schemes.
\item
A {\em line\/} on a smooth cubic surface
$p \colon C \to S$
is a
\mbox{$\bP^1$-bundle}
$l \subset C$
over~$S$
such that, for every
$x \in S$,
one has
$\deg_{\calO(1)} l_x = 1$.
\item
A {\em family of marked cubic surfaces\/} over a base
scheme~$S$
or simply a {\em marked cubic surface\/}
over~$S$
is a cubic surface
$p\colon C \to S$
together with a sequence
$(l_1, \ldots, l_6)$
of six mutually disjoint~lines. The~sequence
$(l_1, \ldots, l_6)$
itself will be called a {\em marking\/}
on\/~$C$.
\end{iii}
\end{defs}

\begin{rems}\smallskip
\begin{iii}
\item
The
\mbox{$\bP^3$-bundle}~$\bP(\calE)$
is not part of the structure of a cubic surface over a base~scheme. Nevertheless,~at least for
$p$~smooth,
we have
$\calO(1) |_C = (\Omega^{\wedge 2}_{C/S})^\vee \!\otimes\! \calL$
for some invertible sheaf
$\calL$
on~$S$.
Thus,~the class of
$\calO(1) |_C$
in
$\Pic(C) / p^*\Pic(S)$
is completely determined by the~datum.
\item
A~marked cubic surface is automatically smooth, according to our~definition. All~its 27~lines are defined
over~$S$.
They~may be labelled as
$l_1, \ldots, l_6$,
$l'_1, \ldots, l'_6$,
$l''_{12}, l''_{13}, \ldots, l''_{56}$,
cf.~\cite[Theorem~V.4.9]{Ha}.
\item
It~is known since the days of A.\,Cayley that there are exactly
$51\,840$
possible markings for a smooth cubic surface with all 27~lines defined over the~base. They~are acted upon, in a transitive manner, by a group of that order, which is isomorphic to the Weyl group
$W(E_6)$~\cite[Theorem~23.9]{Ma}.
\end{iii}
\end{rems}

\begin{conv}\smallskip
In~this article, we will identify
$W(E_6)$
with the permutation group acting on the 27
labels~$l_1, \ldots, l_6, l'_1, \ldots, l'_6, l''_{12}, l''_{13}, \ldots, l''_{56}$.
\end{conv}

\begin{prop}\smallskip
Let\/~$K$
be a~field. Then~there exists a fine moduli
scheme\/~$\smash{\wcalM}$
of marked cubic surfaces
over\/~$K$.
I.e.,~the~functor
\begin{eqnarray*}
F \colon \{K\text{\rm -schemes}\} & \longrightarrow & \{\text{\rm sets}\} \, ,\\
S & \mapsto & \{\text{\rm marked cubic surfaces over~} S\} / \!\sim
\end{eqnarray*}
is representable by a\/
\mbox{$K$-scheme}~$\smash{\wcalM}$.\medskip

\noindent
{\bf Sketch of proof.}
{\em
Let~$\calU \subset (\bP^2)^6$
be the open subscheme parametrizing all ordered
\mbox{$6$-tuples}
of points
on~$\bP^2$
that are in general position. I.e.,~no three lie on a line and not all six lie on a~conic.
$\calU$~is
acted upon, in an obvious manner, by the algebraic
group~$\PGL_3$.

The Hilbert-Mumford numerical criterion~\cite[Theorem~2.1]{MFK} immediately implies that every point
$p \in \calU$
is
\mbox{$\PGL_3$-stable}.
In~fact, the nonstable points
on~$(\bP^2)^6$
are those corresponding to configurations such that there are at least four points on a line~\cite[Definition 3.7/Proposition~3.4]{MFK}. Hence,~the quotient
scheme
$\calU\!/\PGL_3$~exists.

It~is well known that
$\calU\!/\PGL_3$
is the desired fine moduli~scheme. A~formal proof follows the lines of the proof of \cite[Theorem~IV.13]{Be}, with the base field replaced by an arbitrary base~scheme.
}
\eop
\end{prop}

\begin{rems}\smallskip
\begin{iii}
\item
$\smash{\wcalM}$~is
a quasi-projective fourfold. In~fact, such quotients are quasi-projective in much more generality~\cite[Theorem 1.10.ii]{MFK}.
\item
By~functoriality,
$\smash{\wcalM}$
is acted upon
by~$W(E_6)$.
More~precisely, every
$g \in G$
defines a permutation of the 27~labels. For~every base
scheme~$S$,
this defines a map
$T_g(S) \colon F(S) \to F(S)$,
which is natural
in~$S$.
By~Yoneda's lemma, that is equivalent to giving a
morphism~$\smash{T_g\colon \wcalM \to \wcalM}$.
Clearly,~$T_{gg'} = T_g T_{g'}$
for~$g, g' \in W(E_6)$
and
$T_e = \id$
for~$e \in W(E_6)$
the neutral~element.

The~operation
of~$W(E_6)$
is not free, as cubic surfaces may have~automorphisms. It~is, however, free on a non-empty Zariski open subset
of~$\smash{\wcalM}$.
\end{iii}
\end{rems}

\begin{rems}[{\rm A naive embedding}{}]\smallskip
\label{naiv}
\begin{iii}
\item
To~give a
\mbox{$\overline{K}$-rational}
point~$p$
on the variety
$\calU$
is equivalent to giving a sequence of six points
$p_1, \ldots, p_6 \in \bP^2(\overline{K})$
in general~position. A~standard result from projective geometry states that there is a unique
$\gamma \in \PGL_3(\overline{K})$
mapping
$(p_1, p_2, p_3, p_4)$
to the standard basis
$((1:0:0), (0:1:0), (0:0:1), (1:1:1))$.
The~\mbox{$\overline{K}$-rational}
points
on~$\smash{\wcalM}$
may thus be represented by
$3 \times 6$-matrices
of the~form
$$
\left(
\begin{array}{cccccc}
1\, & 0\, & 0\, & 1 & w & y \\
0\, & 1\, & 0\, & 1 & x & z \\
0\, & 0\, & 1\, & 1 & 1 & 1
\end{array}
\right) .
$$
Observe~that vanishing of the third coordinate of
$p_5$
would mean that
$p_1, p_2$,
and~$p_5$
were collinear, and similarly
for~$p_6$.
Hence,~we actually have an open embedding
$\smash{\wcalM \hookrightarrow \bA^4}$.
\item
In particular, one sees that
$\smash{\wcalM}$
is a smooth, affine~scheme. Indeed,~the image of the naive embedding of
$\smash{\wcalM}$
in~$\bA^4$
is the complement of a~divisor.
\end{iii}
\end{rems}

\begin{rem}[{\rm Cayley's compactification}{}]\smallskip
\label{hist1}
The moduli scheme
$\smash{\wcalM}$
of marked cubic surfaces has its origins in the middle of the 19th~century. In~principle, it appears in the article~\cite{Ca} of Arthur~Cayley. Cayley's~approach was as~follows.

Every~smooth cubic surface over an algebraically closed field has 45 tritangent planes meeting the surface in three~lines. Through~each line there are five tritangent~planes. This~leads to a total of 135 cross ratios, which are invariants of the cubic surface, as soon as a marking is fixed on the~lines.

It~turns out that only 45 of these cross ratios are essentially different, due to constraints within the cubic~surfaces. Furthermore,~they provide an embedding
$\smash{\wcalM \hookrightarrow (\bP^1)^{45}}$.
The~image is Cayley's ``cross ratio variety''. For~a more recent treatment of this compactification, we refer the reader to I.~Naruki~\cite{Na}.
\end{rem}

\section{Coble's compactification. The gamma variety.}

\paragraph{\it Coble's irrational invariants.}

\begin{ttt}\smallskip
An~advantage of the algebraic group
$\SL_3$
over the
group~$\PGL_3$
is that its operation on
$\bP^2$
is~linear. This~means that
$\SL_3$
operates naturally on
$\calO(n)$,
and hence on
$\Gamma(\bP^2, \calO(n))$,
for
every~$n$.
It~is well known that there is no
\mbox{$\PGL_3$-linearization}
for
$\calO_{\bP^2}(1)$~\cite[Chapter~1, \S3]{MFK}.

There~is, however, the canonical isogeny
$\SL_3 \twoheadrightarrow \PGL_3$,
the kernel of which consists of the multiples of the identity matrix by the third roots of~unity. These~matrices clearly operate trivially
on~$\calO(3)$.
Thus,~there is a canonical
\mbox{$\PGL_3$-linearization}
for
$\calO(3)$,
which is compatible with the
\mbox{$\SL_3$-linearization}, cf.~\cite[Chapter~3, \S1]{MFK}.

We~may also speak
of~\mbox{$\SL_3$-{\em invariant\/}}
sections of the outer tensor products
$\calO(n_1) \boxtimes \ldots \boxtimes \calO(n_6)$
on
$(\bP^2)^6$
for~$(n_1, \ldots, n_6) \in \bbZ^6$.
If~$3 \,|\, n_1, \ldots, n_6$
then
$\PGL_3$
operates, too, and the
\mbox{$\PGL_3$-invariant}
sections are the same as the
\mbox{$\SL_3$-invariant}~ones.
\end{ttt}

\begin{ttt}\smallskip
For~example, for
$1 \leq i_1 < i_2 < i_3 \leq 6$,
the corresponding minor
$$m_{\{i_1,i_2,i_3\}} := \det
\left(
\begin{array}{ccc}
x_{i_1,0} & x_{i_1,1} & x_{i_1,2} \\
x_{i_2,0} & x_{i_2,1} & x_{i_2,2} \\
x_{i_3,0} & x_{i_3,1} & x_{i_3,2}
\end{array}
\right)
$$
of the
$6 \times 3$-matrix
$$
\left(
\begin{array}{ccc}
x_{1,0} & x_{1,1} & x_{1,2} \\
x_{2,0} & x_{2,1} & x_{2,2} \\
        & \ldots  &         \\
x_{6,0} & x_{6,1} & x_{6,2}
\end{array}
\right)
$$
defines an invariant section
of~$\calO(n_1) \boxtimes \ldots \boxtimes \calO(n_6)$
for
$$n_i :=
\left\{
\begin{array}{cc}
1 & {\rm for~} i \in \{i_1,i_2,i_3\}\, , \\
0 & {\rm for~} i \not\in \{i_1,i_2,i_3\}\, .
\end{array}
\right.$$

Further,
$$
{
\arraycolsep2.5pt
d_2 := \det\!
\left(\!
\begin{array}{cccccc}
x_{1,0}^2 & x_{1,1}^2 & x_{1,2}^2 & x_{1,0}x_{1,1} & x_{1,0}x_{1,2} & x_{1,1}x_{1,2} \\
x_{2,0}^2 & x_{2,1}^2 & x_{2,2}^2 & x_{2,0}x_{2,1} & x_{2,0}x_{2,2} & x_{2,1}x_{2,2} \\
          &           &           & \ldots         &                &                \\
x_{6,0}^2 & x_{6,1}^2 & x_{6,2}^2 & x_{6,0}x_{6,1} & x_{6,0}x_{6,2} & x_{6,1}x_{6,2}
\end{array}
\!\right)
\!\!\in\! \Gamma((\bP^2)^6, \calO(2) \!\boxtimes\! \ldots \!\boxtimes\! \calO(2))
}
$$
is
\mbox{$\SL_3$-invariant},~too.

A.\,Coble~\cite[formulas (16) and~(18)]{Co3} now defines 40
\mbox{$\SL_3$-invariant},
and hence
\mbox{$\PGL_3$-invariant},
sections
$\gamma_. \in \Gamma((\bP^2)^6, \calO(3) \boxtimes \ldots \boxtimes \calO(3))$.
\end{ttt}

\begin{defio}[{\rm Coble}{}]\smallskip
For~$\{i_1, \ldots, i_6\} = \{1, \ldots, 6\}$,
consider
\begin{eqnarray*}
\gamma_{(i_1i_2i_3)(i_4i_5i_6)} & := & m_{\{i_1,i_2,i_3\}} m_{\{i_4,i_5,i_6\}} \,d_2 \hspace{2.5cm} {\rm and} \\
\gamma_{(i_1i_2)(i_3i_4)(i_5i_6)} & := & m_{\{i_1,i_3,i_4\}} m_{\{i_2,i_3,i_4\}} m_{\{i_3,i_5,i_6\}} m_{\{i_4,i_5,i_6\}} m_{\{i_5,i_1,i_2\}} m_{\{i_6,i_1,i_2\}} \, .
\end{eqnarray*}
Following~the original work, we will call these 40~sections the {\em irrational~invariants}.
\end{defio}

\begin{rems}\smallskip
\begin{iii}
\item
Here,~the combinatorial structure is as~follows. Within the parentheses, the indices may be arbitrarily permuted without changing the~symbol. Further,~in the symbols
$\gamma_{(i_1i_2i_3)(i_4i_5i_6)}$,
the two triples may be~interchanged. However,~in the symbols
$\gamma_{(i_1i_2)(i_3i_4)(i_5i_6)}$,
the three pairs may be permuted only~cyclically. Thus,~altogether, there are ten invariants of the first type and 30 invariants of the second~type.
\item
The~20 minors
$m_{\{i_1,i_2,i_3\}}$
and the
invariant~$d_2$
vanish only when the underlying six points
$(x_1, \ldots, x_6)$
are not in general~position.
Hence,~on~$\calU$,
Coble's 40 sections have no~zeroes.
\item
One has the beautiful relation
$$d_2 = - \det
\left(
\begin{array}{cc}
m_{\{1,3,4\}} m_{\{1,5,6\}} & m_{\{1,3,5\}} m_{\{1,4,6\}} \\
m_{\{2,3,4\}} m_{\{2,5,6\}} & m_{\{2,3,5\}} m_{\{2,4,6\}}
\end{array}
\right) ,
$$
cf.~\cite[(47)]{Co1} or \cite[formula~(4.18)]{Hu}.
\end{iii}
\end{rems}

\begin{cau}\smallskip
We~have 40 sections
$\gamma_. \in \Gamma((\bP^2)^6, \calO(3) \boxtimes \ldots \boxtimes \calO(3))^{\PGL_3}$
and a machine calculation shows
$\dim \Gamma((\bP^2)^6, \calO(3) \boxtimes \ldots \boxtimes \calO(3))^{\PGL_3} = 40$.

It~is, however, long known~\cite[(24)]{Co3} that the 40 sections
$\gamma_.$
span only a subvector space of dimension~ten. The~mere fact that there is such a gap is quite~obvious. In~fact, for
$(p_1, \ldots, p_6) \in (\bP^2)^6$
such that
$p_1, \ldots p_4$
are distinct points on a
line~$l$
and
$p_5, p_6 \not\in l$,
we have
$m_{\{1,2,5\}}^3 m_{\{3,4,6\}}^3 \neq 0$
but all
$\gamma_.$~vanish.

In~particular,~the irrational invariants
$\gamma_.$
do not generate the invariant ring
$$\bigoplus_{d \geq 0} \Gamma((\bP^2)^6, \calO(3d) \boxtimes \ldots \boxtimes \calO(3d))^{\PGL_3}$$
and do not define an embedding of the categorical quotient
$((\bP^2)^6)^{\text{\rm ss}} / \PGL_3$
\cite[Definition~0.5]{MFK}
into~$\bP^{39}$.
Observe,~however, Theorem~\ref{Cobleemb}~below.
\end{cau}

\begin{nota}\smallskip
The\/~\mbox{$\PGL_3$-invariant}
local sections
of\/~$\calO(3) \boxtimes \ldots \boxtimes \calO(3)$
form an invertible sheaf
on\/~$\smash{\wcalM = \calU\!/\PGL_3}$,
which we will denote
by\/~$\calL$.
\end{nota}

\begin{theo}\smallskip
\label{Cobleemb}
\begin{abc}
\item
The~invertible
sheaf\/~$\calL$
on\/~$\smash{\wcalM}$
is very~ample.
\item
The~40 irrational invariants\/
$\smash{\gamma_. \in \Gamma(\wcalM, \calL)}$
define a projective embedding\/
$\smash{\gamma\colon \wcalM \hookrightarrow \bP_K^{39}}$.
\item
The~Zariski closure\/
$\smash{\widetilde{M}}$
of the image
of\/~$\gamma$
is contained in a nine-dimensional linear~subspace.
\item
As~a subvariety of this\/
$\bP^9$,
$\smash{\widetilde{M}}$
has the properties~below.
\begin{iii}
\item
The~image
of\/~$\smash{\widetilde{M}}$
under the\/
$2$-uple
Veronese embedding\/
$\bP^9 \hookrightarrow \bP^{54}$
is not contained in any proper linear~subspace.
\item
The~image
of\/~$\smash{\widetilde{M}}$
under the\/
$3$-uple
Veronese embedding\/
$\bP^9 \hookrightarrow \bP^{219}$
is contained in a linear subspace of
dimension\/~$189$.
\item
$\smash{\widetilde{M}}$
is the intersection of 30 cubic~hypersurfaces.
\end{iii}
\end{abc}\medskip

\noindent
{\bf Proof.}
{\em
Assertion~b) is \cite[Corollary~5.9]{CGL}. The~proof given there is based on the considerations of I.~Naruki~\cite{Na}. a)~is clearly implied by~b).\smallskip

\noindent
c)~follows from the fact that the vector space
$\langle \gamma_. \rangle$
spanned by the 40 irrational invariants
$\gamma_.$
is ten-dimensional.\smallskip

\noindent
d.i) and~ii) As~is easily checked by computer, the purely quadratic expressions in the
$\gamma_.$
form a 55-dimensional vector space, while the purely cubic expressions form a vector space of dimension~190.\smallskip

\noindent
iii)
By~ii),
$\smash{\widetilde{M}}$
is contained in the intersection of 30 cubic hypersurfaces
in~$\bP^9$.
This~intersection is reported by {\tt magma} as being reduced and irreducible of dimension~four.
}
\eop
\end{theo}

\begin{rems}\smallskip
\begin{iii}
\item
As~$\smash{\gamma\colon \wcalM \to \bP_K^{39}}$
is a map that is given completely explicitly, one might try to use computer algebra to prove it is an~embedding. This~indeed works, at least when one organizes the computation in a slightly deliberate~way.

It~turns out that the composition
$\widetilde\gamma$
of~$\gamma$
with the linear projection to
the~$\bP^9$,
formed by the ten invariants of
type~$\gamma_{(i_1i_2i_3)(i_4i_5i_6)}$
is already defined everywhere and separates tangent vectors. It~is a
$2\!:\!1$-morphism
identifying
$(w,x,y,z)$
with
$(w',x',y',z')$~for
$$
\scriptstyle
w' := \frac{(wz - xy)(z-1)}{(x-z)(y-z)},~
x' := \frac{(wz - xy)(y-1)}{(w-y)(y-z)},~
y' := \frac{(wz - xy)(x-1)}{(w-x)(x-z)},
\text{ and }
z' := \frac{(wz - xy)(w-1)}{(w-x)(w-y)} \, .
$$
A~Gr\"obner base calculation in four variables readily shows that these two points never have the same image
under~$\gamma$.
\item
The partner point
$(w',x',y',z')$
corresponds to the same cubic surface
as
$(w,x,y,z)$,
but with the flipped~marking.
I.e.,~$l_i$
is replaced by
$l'_i$
and vice~versa. This~is seen by a short calculation from \cite[Table~(2)]{Co3}, cf.~\cite[p.~196]{Co1}.
\end{iii}
\end{rems}

\begin{defs}\smallskip
\begin{iii}
\item
We~will call
$\smash{\gamma\colon \wcalM \hookrightarrow \bP_K^{39}}$
{\em Coble's gamma~map}.
\item
The~variety
$\smash{\widetilde{M}}$,
given as the Zariski closure of the image
of~$\gamma$
will be called {\em Coble's gamma~variety.}
\end{iii}
\end{defs}

\begin{rems}\smallskip
\begin{iii}
\item
The~fact that the vector space
$\langle \gamma_. \rangle$
is only of dimension ten is, of course, easily checked by computer, as~well.
\item
The assertions c) and~d) are due to A.~Coble~himself. For~d), we advise the reader to compare the result \cite[Theorem~6.4]{AF} of D.~Allcock and E.~Frei\-tag.

Coble's~original proof for~c) works as follows~\cite[(24)]{Co3}. One~may write down~\cite[page~343]{Co3} five four-term linear relations, the
\mbox{$S_6$-orbits}
of which yield a total of 270~relations. These~relations form a single orbit
under~$W(E_6)$
and generate the 30-dimensional space of all linear~relations.

In~order to show that the dimension is not lower than ten, Coble has to use the moduli~interpretation. He~verifies that there are enough cubic surfaces in hexahedral~form.
\item
The~cubic relations are in fact more elementary than the linear~ones. For~example, one~has
$$\gamma_{(12)(34)(56)} \gamma_{(23)(45)(16)} \gamma_{(14)(36)(25)} = \gamma_{(12)(36)(45)} \gamma_{(34)(25)(16)} \gamma_{(56)(14)(23)} \, .$$
To~see this, look at the left hand side~first. The~nine pairs of numbers in
$\{1, \ldots, 6\}$
that are used, are exactly those with an odd~difference. Thus,~when writing, according to the very definition, the left side as a product of 18~minors,
$m_{1,3,5}$
and
$m_{2,4,6}$
can not~appear. It~turns out that each of the other minors occurs exactly~once. As~the same is true for the right hand side, the equality becomes~evident.

We~remark that this relation is not a consequence of the linear~ones. I.e., it does not become trivial when restricted
to~$\bP^9$.
Its~orbit under
$W(E_6)$
must generate the 30-dimensional space of all cubic~relations. Indeed,~that is an irreducible representation, as we will show in the next~subsection.
\item
In~particular, the gamma variety
$\widetilde{M}$
is clearly not a complete~intersection. Nevertheless,~the following of its numerical invariants may be~computed.
\end{iii}
\end{rems}

\begin{lem}\smallskip
\label{CastMum}
\begin{iii}
\item
The Hilbert series
of\/~$\smash{\widetilde{M}}$
is\/
$\smash{\frac{1+5T+15T^2+5T^3+T^4}{(1-T)^5}}$.
\item
In~particular, the Hilbert polynomial
of\/~$\smash{\widetilde{M}}$
is\/
$\smash{\frac98 T^4 + \frac94 T^3 + \frac{27}8 T^2 + \frac94 T + 1}$.
Further,~the Hilbert polynomial agrees with the Hilbert function in all
degrees\/~$\geq \!0$.
\item
$\smash{\widetilde{M}}$~is
a projective variety of
degree\/~$27$.
\item
The~Castelnuovo-Mumford regularity
of\/~$\smash{\widetilde{M}}$
is equal
to\/~$4$
and that of the ideal sheaf\/
$\calI_{\widetilde{M}} \subset \calO_{\bP^{39}}$
is equal
to\/~$5$.
\end{iii}\medskip

\noindent
{\bf Proof.}
{\em
i)~follows from a Gr\"obner base~calculation. ii)~and iii) are immediate con\-se\-quences of~i).\smallskip

\noindent
iv)
By~\cite[p.\,219]{DE}, it is pure linear algebra to compute the Castelnuovo-Mumford regularity of a coherent
$\calO_{\bP^N}$-module.
We~used the implementation in~{\tt magma}.
}
\eop
\end{lem}\medskip

\paragraph{\it The operation of\/~$W(E_6)$.}

\begin{ttt}\smallskip
\label{we6lin}
It~is an important feature of Coble's (as well as Cayley's) compactifications that they explicitly linearize the operation
of~$W(E_6)$.
More~precisely,\medskip

\noindent
{\bf Lemma.}
{\em
There~exists a\/
\mbox{$W(E_6)$-linearization}
of\/~$\smash{\calL \in \Pic(\wcalM)}$
such that

\begin{iii}
\item
the 80~sections\/
$\smash{\pm\gamma_. \in \Gamma(\wcalM, \calL)}$
form a\/
\mbox{$W(E_6)$-invariant}~set.
\item
The~corresponding permutation representation\/
$\Pi\colon W(E_6) \hookrightarrow S_{80}$
is~transitive. It~has a system of 40 blocks given by the
pairs\/~$\{\gamma, -\gamma\}$.
\item
The~permutation representation\/
$W(E_6) \hookrightarrow S_{40}$
on the 40~blocks is the same as that on decompositions of the 27~lines into three pairs of Steiner~trihedra.
\end{iii}\smallskip

\noindent
{\bf Proof.}
{\em
i)
(Cf.~\cite[Section~2]{CGL})
As~$W(E_6)$
is a discrete group, the general concept of a linearization of an invertible sheaf~\cite[Definition~1.6]{MFK} breaks down to a system of compatible isomorphisms
$\smash{i_g\colon T_g^* \calL \stackrel{\cong}{\longrightarrow} \calL}$
for
$\smash{T_g\colon \wcalM \to \wcalM}$
the operation
of~$g$.

For~$g \in S_6 \subset W(E_6)$,
there is an obvious such~isomorphism.
Indeed,~$g$
permutes the six labels
$l_1, \ldots, l_6$
and, accordingly, the six blow-up
points~$p_1, \ldots, p_6$.
Simply~permute the six factors
of~$\calO(3) \boxtimes \ldots \boxtimes \calO(3)$
as described
by~$g$.
Assertion~i) is clear for these~elements.

Further,~$W(E_6)$
is generated
by~$S_6$
and just one additional element, the quadratic
transformation~$I_{123}$
with centre in
$p_1$,
$p_2$,
and
$p_3$~\cite[Example~V.4.2.3]{Ha}.
In~the coordinates described in Remark~\ref{naiv}.i), this map is given
by~$\smash{(w,x,y,z) \mapsto (\frac1w, \frac1x, \frac1y, \frac1z)}$.

One~may now list explicit formulas for the 40 irrational invariants
$\gamma_.$
in terms of these~coordinates. Each~of these sections actually defines a global trivialization
of~$\calL$.
Plugging~in the provision
$\smash{(w,x,y,z) \mapsto (\frac1w, \frac1x, \frac1y, \frac1z)}$
in a naive way, yields an isomorphism
$\smash{i'_{I_{123}}\colon T_{I_{123}}^* \calL \stackrel{\cong}{\longrightarrow} \calL}$.
It~turns out that,
under~$i'_{I_{123}}$,
the 40 sections
$\gamma_.$
are permuted up to signs and a common scaling factor
of~$\smash{\frac1{w^2x^2y^2z^2}}$.
Thus,~let us take
$\smash{i_{I_{123}} := w^2x^2y^2z^2 \!\cdot i'_{I_{123}}}$
as the actual~definition.

This~uniquely determines
$i_g$
for
every~$g \in W(E_6)$.
One~may check that
$\{i_g\}_{g \in W(E_6)}$
is a well-defined linearization
of~$\calL$.
Assertion~i) is then~clear.\smallskip

\noindent
ii)
We~checked the first assertion in~{\tt magma}. The~second statement is~obvious.\smallskip

\noindent
iii)
Note~that, in the blown-up model, the 40 irrational invariants have exactly the same combinatorial structure as the 40 decompositions, cf.~\cite[3.7]{EJ5}.
}}
\eop
\end{ttt}

\begin{rems}\smallskip
\begin{iii}
\item
The~permutation representation
$\Pi$
has no other nontrivial block structures.
\item
The~restriction
of~$\Pi$
to the index-two subgroup
$D^1 W(E_6) \subset W(E_6)$,
which is the simple group of
order~$25\,920$,
is still~transitive. Neither~does it have more block~structures.
\item
Lemma~\ref{we6lin}.i) suggests that it might have technical advantages to consider the embedding
$\smash{\gamma'\colon \wcalM \hookrightarrow \bP^{79}}$,
linearly equivalent to the gamma
map~$\gamma$,
which is defined by the 80
sections~$\pm\gamma_.$.
To~a certain extent, this is indeed the case, cf.~Remarks~\ref{expldesc}~below.
\end{iii}
\end{rems}

\begin{remso}[{\rm Representations of~$W(E_6)$}{}]\smallskip
\begin{iii}
\item
The~dimensions of the irreducible complex representations
of~$W(E_6)$
are
$1$, $1$, $6$, $6$, $10$, $15$, $15$, $15$, $15$, $20$, $20$, $20$, $24$, $24$, $30$, $30$, $60$, $60$, $60$, $64$, $64$, $80$, $81$, $81$,
and~$90$.
\item
The~\mbox{$W(E_6)$-representation}
on the vector space
$V := \langle \gamma_. \rangle \cong \Gamma(\bP^9, \calO(1))$
of dimension ten is~irreducible.
\item
The~\mbox{$W(E_6)$-representation}
on the 220-dimensional vector space
$\Gamma(\bP^9, \calO(3))$
decomposes into two copies of the ten-dimensional, two copies of a 30-dimensional, two copies of the other 30-dimensional, and one copy of the 80-dimensional irreducible representations~\cite[Theorem~3.2.2]{DK}. This~already implies that the 30-dimensional sub-representation of cubic relations among the
$\gamma_.$
is~irreducible.
\end{iii}
\end{remso}

\begin{rems}\smallskip
\label{hist2}
\begin{abc}
\item
The~embedding of the moduli scheme of marked cubic surfaces
into~$\bP^9$,
originally due to A.\,B.~Coble, was studied recently by D.\,Allcock and E.\,Freitag~\cite{AF}, as well as B.\,van Geemen~\cite{Ge}. Their~approaches were rather different from Coble's. For~example, van Geemen actually constructs an embedding of the cross ratio variety, instead
of~$\calU\!/\PGL_3$,
into~$\bP^9$.
He~obtains the 30~cubic relations in~\cite[7.9]{Ge}.
\item
A~short summary of Coble's approach may be found in I.~Dolgachev's book on classical algebraic geometry~\cite[Remark~9.4.20]{Do}.
\end{abc}
\end{rems}

\section{The moduli scheme of un-marked cubic surfaces}

\begin{ttt}\smallskip
The~quotient
$\smash{\wcalM / W(E_6) =: \calM}$
is the coarse moduli scheme of smooth cubic~surfaces. The~reader might consult~\cite[Appendix by E.\,~Looijenga]{Na} for more details on this~quotient. As~cubic surfaces may have automorphisms, a fine moduli scheme cannot~exist.
\end{ttt}

\begin{ttt}\smallskip
The~moduli scheme of smooth cubic~surfaces may as well be constructed directly as the quotient
$\calV / \PGL_4$,
for
$\calV \subset \bP(\Sym^3(K^4)^*) \cong \bP^{19}$
the open subscheme parametrizing smooth cubic~surfaces. In~fact, by~\cite[1.14]{Mu2}, every smooth cubic surface corresponds to a
\mbox{$\PGL_4$-stable}
point
on~$\bP^{19}$.

The~\mbox{$\PGL_4$-invariants}
have been determined by A.\,Clebsch~\cite[sections 4 and~5]{Cl} as early as~1861. In~today's language, Clebsch's result is that there is an open embedding
$\Cl\colon \calV / \PGL_4 \cong \calM \hookrightarrow \bP(1,2,3,4,5)$
into a weighted projective space~\cite[formula~(9.57)]{Do}.
\end{ttt}

\begin{defi}\smallskip
\begin{iii}
\item
The~homogeneous coordinates on
$\bP(1,2,3,4,5)$
will be denoted, in this order, by
$A$,
$B$,
$C$,
$D$,
and~$E$.
\item
Thus,~given a smooth cubic surface over a
field~$K$,
there is the corresponding
\mbox{$K$-rational}
point
on~$\bP(1,2,3,4,5)$.
Its~homogeneous coordinates form a vector
$[A, \ldots, E]$,
which is unique up to weighted scaling, for the weight
vector~$(1, \ldots, 5)$.
We~will speak of {\em Clebsch's invariant vector\/} or simply {\em Clebsch's invariants\/} of the cubic~surface.
\end{iii}
\end{defi}

\begin{ex}\smallskip
Consider the {\em pentahedral\/} family
$\calC \to \bP^4/S_5$
of cubic surfaces, given by
\begin{eqnarray}
\label{pent}
a_0 X_0^3 + a_1 X_1^3 + a_2 X_2^3 + a_3 X_3^3 + a_4 X_4^3 & = & 0 \,, \\
\phantom{a_0}\hsmash{X_0}\phantom{X_0^3} + \phantom{a_1}\hsmash{X_1}\phantom{X_1^3} + \phantom{a_2}\hsmash{X_2}\phantom{X_2^3} + \phantom{a_3}\hsmash{X_3}\phantom{X_3^3} + \phantom{a_4}\hsmash{X_4}\phantom{X_4^3} & = & 0 \nonumber
\end{eqnarray}
over~$\bP^4/S_5 \cong \bP(1,2,3,4,5)$.
We~will use the elementary symmetric functions
$\sigma_1, \ldots, \sigma_5$
in~$a_0, \ldots, a_4$
as natural homogeneous coordinates
on~$\bP^4/S_5$.

Let~us restrict our considerations to the open subset
$\calP \subset \bP^4/S_5$
representing smooth cubic surfaces having a {\em proper pentahedron}. The~latter condition is equivalent
to~$\sigma_5 \neq 0$.

Then,~for
$t\colon \calP \longrightarrow \calM$
the classifying morphism, the composition
$\smash{\Cl \!\circ\! t\colon \calP \to \calM \hookrightarrow \bP(1,2,3,4,5)}$
is given by the
\mbox{$S_5$-invariant}
sections
\begin{equation}
\label{Clebsch}
I_8 := \sigma_4^2 - 4\sigma_3\sigma_5, \quad I_{16} := \sigma_1\sigma_5^3, \quad I_{24} := \sigma_4\sigma_5^4, \quad I_{32} := \sigma_2\sigma_5^6, \quad I_{40} := \sigma_5^8
\end{equation}
of~$\calO(8)$,
$\calO(16)$,
$\calO(24)$,
$\calO(32)$,
and
$\calO(40)$,
respectively. See~\cite[formula~(9.59)]{Do} or \cite[paragraph~543]{Sa}. In~other words
$(\Cl \!\circ\! t)^{-1} (A) = I_8$,
\ldots,
$(\Cl \!\circ\! t)^{-1} (E) = I_{40}$.
\end{ex}

\begin{lem}\smallskip
The classifying morphism\/
$t\colon \calP \to \calM$
is an open embedding.\smallskip

\noindent
{\bf Proof.}
{\em
It~will suffice to show that
$\Cl \!\circ\! t \colon \calP \to \bP(1,2,3,4,5)$
is an open~embedding. For~this, we first observe that
$\Cl \!\circ\! t$
is~birational. Indeed,~the two function fields are
\begin{eqnarray*}
K(\bP(1,2,3,4,5)) & = & K(B/A^2, C/A^3, D/A^4, E/A^5) \\
 & = & K(A^2/B, A^3/C, A^4/D, A^5/E)
\end{eqnarray*}
and
$K(\calP) = K(\sigma_2/\sigma_1^2, \sigma_3/\sigma_1^3, \sigma_4/\sigma_1^4, \sigma_5/\sigma_1^5)$.
Both~are of transcendence degree four
over~$K$.

Consider the finitely generated
\mbox{$K$-algebra}
$$R := K[\frac{A^2}B, \!\frac{A^3}C, \!\frac{A^4}D\!, \!\frac{A^5}E\!, \!\frac{D}{B^2}, \!\frac{C^2-AE}{4B^3}, \!\frac{CE}{B^4}, \!\frac{E^2}{B^5}] \, ,$$
which is a subdomain
of~$K(\bP(1,2,3,4,5))$.
The~formulas~(\ref{Clebsch}) together with
\begin{align*}
\textstyle
(\Cl \!\circ\! t)^{-1} (\frac{D}{B^2}) = \sigma_2/\sigma_1^2& \, , \quad
\textstyle
(\Cl \!\circ\! t)^{-1} (\frac{C^2-AE}{4B^3}) = \sigma_3/\sigma_1^3 \, , \\
\textstyle
(\Cl \!\circ\! t)^{-1} (\frac{CE}{B^4}) = \sigma_4/\sigma_1^4& \, , \quad
\textstyle
\hspace{0.7cm} (\Cl \!\circ\! t)^{-1} (\frac{E^2}{B^5}) = \sigma_5/\sigma_1^5 \, ,
\end{align*}
immediately define a
\mbox{$K$-algebra}
homomorphism
$\iota\colon R \to K(\calP)$.
For~$\frakp := \ker\iota$,
we have a homomorphism
$\Q(R/\frakp) \hookrightarrow K(\calP)$
of~fields.

As~$\sigma_2/\sigma_1^2$,
$\sigma_3/\sigma_1^3$,
$\sigma_4/\sigma_1^4$,
and~$\sigma_5/\sigma_1^5$
are in the image, we see that
$(\Cl \!\circ\! t)^{-1}$
actually~defines an isomorphism
$\Q(R/\frakp) \cong K(\calP)$.
In~particular,
$\Q(R/\frakp)$
is of transcendence degree four and,
consequently,~$\frakp = (0)$.
As~$\Q(R) = K(\bP(1,2,3,4,5))$,
the claim~follows.

Furthermore,~$\Cl \!\circ\! t$
is a quasi-finite~morphism. In~fact, this may be tested on closed points and after base extension to the algebraic
closure~$\overline{K}$.
Thus,~let
$p = (A, \ldots, E) \in \bP(1,2,3,4,5)(\overline{K})$
be a geometric~point.
If~$E = 0$
then
$(\Cl \!\circ\! t)^{-1}(p) = \emptyset$.
Otherwise,~there are eight solutions of
$\sigma_5^8 = E$
and, for each choice,
$\sigma_1, \ldots, \sigma_4$
may be computed~directly.

Finally,~$\bP(1,2,3,4,5)$
is a toric variety~\cite[section~2.2, page~35]{Fu} and hence a normal scheme~\cite[section~2.1, page~29]{Fu}. Therefore~the assertion is implied by \cite[Corollaire~(4.4.9)]{EGAIII}.
}
\eop
\end{lem}

\begin{rems}\smallskip
\begin{iii}
\item
In~particular, a general cubic surface over a field has a proper~pentahedron, which will usually be defined over a finite extension~field.
\item
Further,~on the open subset
of~$\calM$
representing smooth cubic surfaces with a proper pentahedron,
$\sigma_1, \ldots, \sigma_5$
serve well as~coordinates. It~is highly remarkable that they do not extend properly to the whole
of~$\calM$.
\end{iii}
\end{rems}

\begin{ex}\smallskip
There~are other prominent families of smooth cubic~surfaces. The~most interesting ones are probably the hexahedral~families. Consider
$\calC \to H \subset \bP^5$,
where
$\calC \subset H \times \bP^4$
is given by
\begin{eqnarray*}
\phantom{a_0} X_0^3 + \phantom{a_1} X_1^3 + \phantom{a_2} X_2^3 + \phantom{a_3} X_3^3 + \phantom{a_4} X_4^3 + \phantom{a_5} X_5^3 & = & 0 \,
, \\
\phantom{a_0}\hsmash{X_0}\phantom{X_0^3} + \phantom{a_1}\hsmash{X_1}
\phantom{X_1^3} + \phantom{a_2}\hsmash{X_2}\phantom{X_2^3} + \phantom{a_3}\hsmash{X_3}\phantom{X_3^3} + \phantom{a_4}\hsmash{X_4}\phantom{X_4^3} + \phantom{a_5}\hsmash{X_5}\phantom{X_5^3} & = & 0 \, , \\
a_0 \hsmash{X_0}\phantom{X_0^3} + a_1 \hsmash{X_1}\phantom{X_1^3} + a_2 \hsmash{X_2}\phantom{X_2^3} + a_3 \hsmash{X_3}\phantom{X_3^3} + a_4 \hsmash{X_4}\phantom{X_4^3} + a_5 \hsmash{X_5}\phantom{X_5^3} & = & 0 \, .
\end{eqnarray*}
and
$H \subset \bP^5$
is the hyperplane defined
by~$a_0 + \ldots + a_5 = 0$.
This~is the {\em ordered hexahedral\/} family of cubic~surfaces. Correspondingly,~the base of the {\em unordered hexahedral\/} family is the
quotient~$H/S_6 \cong \bP(2,3,4,5,6)$.

There~are the tautological morphisms
$\smash{\wcalM \stackrel{t_1}{\longrightarrow} H \stackrel{t_2}{\longrightarrow} H/S_6 \stackrel{t_3}{\longrightarrow} \calM}$.
It~is classically known that
$t_1$
is an unramified
$2:1$-covering
and that
$t_3$
is an unramified
\mbox{$36:1$-covering}.
Clearly,~$t_2$
is generically~$720:1$.
\end{ex}

\begin{exo}[{\rm continued}{}]\smallskip
\label{hexa2}
It~seems natural to use the elementary symmetric functions
$\sigma_2, \ldots, \sigma_6$
in the hexahedral coefficients as homogeneous coordinates
on~$H/S_6$.
Then~it is possible, today, to give explicit formulas for the composition
$\smash{\Cl \!\circ\! t_3 \colon H/S_6 \to \calM \hookrightarrow \bP(1,2,3,4,5)}$.

This~means to convert the formulas~(\ref{Clebsch}) for Clebsch's invariants to the hexahedral~form. The~first of these formulas,
\begin{equation}
\label{Sousley}
(\Cl \!\circ\! t_3)^{-1} (A) = 24[4\sigma_2^3 - 3\sigma_3^2 - 16\sigma_2\sigma_4 + 12\sigma_6] \, ,
\end{equation}
was established by C.\,P.\ Sousley~\cite[formula~(17)]{So}, back in~1917. Here,~the
coefficient~$24$
is somewhat conventional, as it depends on the choice of an isomorphism
$(\Cl \!\circ\! t_3)^* \calO(1) \cong \calO(6)$.

Formula~(\ref{Sousley}) agrees with the modern treatment, due to I.\,V.~Dolgachev~\cite[Remark~9.4.19]{Do} as well as with \cite[formula~(B.56)]{Hu}. Other~coefficients were used, however, in Coble's original work~\cite[formula~(9)]{Co3} and to obtain~\cite[formula (4.108)]{Hu}.
\end{exo}

\begin{theo}\smallskip
\label{Formeln}
\begin{iii}
\item
The~canonical morphism\/
$$\smash{\psi \colon \wcalM \stackrel{\pr}{\longrightarrow} \calM \stackrel{\Cl}{\hookrightarrow} \bP(1,2,3,4,5)}$$
allows an extension to\/
$\bP^{39}$
under the gamma~map. More~precisely, there exists a rational map\/
$\widetilde\psi \colon \bP^{39} \ratarrow \bP(1,2,3,4,5)$
such that the following diagram~commutes,
$$
\definemorphism{gleich}\Solid\notip\notip
\definemorphism{incl}\solid\tip\tip
\definemorphism{rat}\dashed\tip\notip
\diagram
\wcalM_{\!\!\phantom g} \rto^\pr \dincl_\gamma & \calM\; \rincl^{\Cl \;\;\;\;\;\;\;\;\;\;\;\;\;\;} & \bP(1,2,3,4,5) \phantom{\, .} \dgleich \\
\bP^{39} \rrrat^{\widetilde\psi \;\;\;\;\;\;\;\;\;\;\;\;} & & \bP(1,2,3,4,5) \, .
\enddiagram
$$
\item
Explicitly,~the rational map\/
$\widetilde\psi \colon \bP^{39} \ratarrow \bP(1,2,3,4,5)$,
defined by the global sections
\begin{eqnarray}
\label{invvsinv}
\bullet & & -6P_2 \in \Gamma(\bP^{39}, \calO(2)) \, , \nonumber \\
\bullet & & \textstyle -24P_4 + \frac{41}{16} P_2^2 \in \Gamma(\bP^{39}, \calO(4)) \, , \nonumber \\
\bullet & & \textstyle \frac{576}{13} P_6 - \frac{396}{13} P_4 P_2 + \frac{29}{13} P_2^3 \in \Gamma(\bP^{39}, \calO(6)) \, , \\
\bullet & & \textstyle -\frac{62208}{1171} P_8 + \frac{54864}{1171} P_6 P_2 + \frac{203616}{1171} P_4^2 - \frac{61287}{1171} P_4 P_2^2 + \frac{13393}{4684} P_2^4 \!\in\! \Gamma(\bP^{39}, \calO(8)) \, , \nonumber \\
\bullet & & \textstyle \frac{41472}{155} P_{10} - \frac{4605984}{36301} P_8 P_2 - \frac{106272}{403} P_6 P_4 + \frac{19990440}{471913} P_6 P_2^2 + \frac{47719206}{471913} P_4^2 P_2 \nonumber \\
& & \textstyle \hspace{4.2cm} {} - \frac{7468023}{471913} P_4 P_2^3 + \frac{10108327}{18876520} P_2^5 \in \Gamma(\bP^{39}, \calO(10)) \, , \nonumber
\end{eqnarray}
satisfies this~condition.
Here,~$P_k$
denotes the sum of the 40
\mbox{$k$-th}~powers.
\item
In~other words, these formulas express Clebsch's invariants\/
$A, \ldots, E$
in terms of Coble's 40 irrational invariants\/~$\gamma_.$.
\end{iii}
\end{theo}

\begin{rems}\smallskip
\begin{iii}
\item
The formula for the first Clebsch invariant is due to A.\,B.~Coble, cf.~\cite[formula~(38)]{Co3} and \cite[formula~(4.108)]{Hu}. It~may be obtained by plugging the formula~\cite[formula (85)]{Co1}, computing hexahedral coefficients out of six blow-up points, into Sousley's formula~(\ref{Sousley}).
\item
Similarly~to \ref{hexa2}, there is a minor ambiguity here, due to the possibility of~scaling. The~coefficient
$(-6)$
in the first formula agrees with Sousley's formula~(\ref{Sousley}).
\end{iii}
\end{rems}

\begin{ttt}\smallskip
{\bf Proof of Theorem~\ref{Formeln}.}
We~will prove this theorem in several~steps.\smallskip

\noindent
{\em First step.}
Results from the literature and preparations.

\noindent
A~considerable part of this result is available from the~literature. First~of all, it is known that the morphism
$\psi := \Cl \!\circ\! \pr$
in the upper row extends to a finite morphism
$\smash{\varphi\colon \widetilde{M} \to \bP(1,2,3,4,5)}$
from the gamma variety~\cite[Proposition~1.3]{CGL}.

Further,
$$\varphi^* \!\calO(1) \cong \calO(2) |_{\widetilde{M}} \, .$$
Indeed,~this follows from the functoriality of the determinant line bundle~\cite[Definition~1.1]{CGL}, together with its calculation for both sides, \cite[Proposition~1.3]{CGL} and \cite[Section~2]{CGL}.

The fact that
$\smash{\varphi\colon \widetilde{M} \to \bP(1,2,3,4,5)}$
is a finite morphism ensures~that
$$\varphi^{-1}\colon \Gamma(\bP(1,2,3,4,5), \calO(i)) \longrightarrow \Gamma(\widetilde{M}, \calO(2i)|_{\widetilde{M}})^{W(E_6)}$$
is a bijection for
each~$i$.
In~particular, every
\mbox{$W(E_6)$-invariant}
even degree homogeneous polynomial expression in Coble's irrational invariants induces an element of Clebsch's invariant~ring.\smallskip

\noindent
{\em Second step.}
Extending the sections
to~$\bP^{39}$.

\noindent
Unfortunately, we need exactly the~opposite. To~ensure this, we claim that the restriction~map
$$\res_i \colon \Gamma(\bP^{39}, \calO(2i))^{W(E_6)} \to \Gamma(\widetilde{M}, \calO(2i)|_{\widetilde{M}})^{W(E_6)}$$
is surjective,
for~$i = 1, \ldots, 5$.

In~view of the bijectivity
of~$\varphi^{-1}$,
it will suffice to verify that
$\dim \im \res_i \geq d_i$~for
$$d_i := \left\{\!
\begin{array}{rr}
\hspace{0.1cm} 1\;\; & \text{\hspace{0.2cm}for } i = 1 \, , \\
\hspace{0.1cm} 2\;\; & \text{\hspace{0.2cm}for } i = 2 \, , \\
\hspace{0.1cm} 3\;\; & \text{\hspace{0.2cm}for } i = 3 \, , \\
\hspace{0.1cm} 5\;\; & \text{\hspace{0.2cm}for } i = 4 \, , \\
\hspace{0.1cm} 7\;\; & \text{\hspace{0.2cm}for } i = 5 \, .
\end{array}
\right.
$$
For~this, let us write down some
\mbox{$W(E_6)$-invariant}
sections
of~$\calO(2i)|_{\widetilde{M}}$
that are contained in the image
of~$\res_i$.
Denote by
$\smash{P_i := \sum_{j=0}^{39} X_j^i}$
the
$i$-th
power~sum.~Then
\begin{align}
\label{base}
&\im(\res_i\colon \Gamma(\bP^{39}, \calO(2i))^{W(E_6)} \to \Gamma(\widetilde{M}, \calO(2i)|_{\widetilde{M}})^{W(E_6)}) \nonumber\\
 \supseteq & \left\{
\begin{array}{lr}
\langle P_2 \rangle & \hspace{0.5cm}\text{for } i = 1 \, , \\
\langle P_4, P_2^2 \rangle & \text{for } i = 2 \, , \\
\langle P_6, P_4P_2, P_2^3 \rangle & \text{for } i = 3 \, , \\
\langle P_8, P_6P_2, P_4^2, P_4P_2^2, P_2^4 \rangle & \text{for } i = 4 \, , \\
\langle P_{10}, P_8P_2, P_6P_4, P_6P_2^2, P_4^2P_2, P_4P_2^3, P_2^5\rangle & \text{for } i = 5 \, .
\end{array}
\right.
\end{align}
It~suffices to verify that, for
$i = 1, \ldots, 5$,
the
$d_i$
global sections given of
$\calO(2i)|_{\widetilde{M}}$
are linearly independent.

These~are simple machine~calculations. Starting~with six
\mbox{$\bbQ$-rational}
points
on~$\bP^2$
in general position, one may compute Coble's irrational invariants and obtains a point
on~$\smash{\widetilde{M} \subset \bP^{39}}$.
Evaluating~the power sums and the expressions listed yields a vector
in~$\bbQ^{d_i}$.
Having~repeated this process
$N$~times,
one ends up with a
\mbox{$d_i \times N$-matrix}
and the task is to show that it is of
rank~$d_i$.
The~calculation may be executed over the rationals or modulo a prime of moderate~size.

The~linear
maps~$\res_i$
are thus surjective and we actually found bases for
$\smash{\Gamma(\widetilde{M}, \calO(2i)|_{\widetilde{M}})^{W(E_6)}}$,
consisting of sections extending to the whole
of~$\bP^{39}$.
This~is enough to prove assertion~i).\smallskip

\noindent
{\em Third step.}
The proofs of ii) and~iii).

\noindent
The rational map
$\varphi$
is defined by five sections
$\smash{s_i \in \Gamma(\widetilde{M}, \calO(2i)|_{\widetilde{M}})^{W(E_6)}}$,
for
$i = 1, \ldots, 5$.
To~explicitly describe an extension
to~$\bP^{39}$
as desired, the actual coefficients
of~$s_1, \ldots, s_5$
in the bases~(\ref{base}) have to be~determined.

This~is, in fact, an interpolation~problem. Starting~with a smooth cubic surface in the blown-up model, one may, as in the second step, directly compute the values of the 40 irrational invariants
$\gamma_.$
and their power~sums. On~the other hand, using the methods described in~\ref{meth1}, Algorithm~\ref{meth2}, and~\ref{meth3}, it is typically possible to compute Clebsch's
invariants~$A, \ldots, E$.
Having~done this for sufficiently many surfaces, the 18~coefficients are fixed up to the appropriate scaling~factors.\smallskip

\noindent
iii)~is only a reformulation of~ii).
\eop
\end{ttt}

\begin{rems}\smallskip
\begin{iii}
\item
We~find it quite noteworthy that the polynomials representing
$s_1, \ldots, s_5$
may actually be chosen to be
$S_{40}$-invariant,
particularly in view of the fact that
$W(E_6)$
is of an enormous index
in~$S_{40}$.
\item
At~least
for~$i \geq 2$,
the full restriction homomorphism
$$\Gamma(\bP^{39}, \calO(2i)) \longrightarrow \Gamma(\widetilde{M}, \calO(2i)|_{\widetilde{M}})$$
is surjective, as may be shown by the usual cohomological~argument. Recall~from Lemma~\ref{CastMum}.iii) that the Castelnuovo-Mumford regularity of
$\calI_{\widetilde{M}}$
is equal
to~$5$.
Since~$1+2i \geq 5$,
this implies
$H^1(\bP^{39}, \calI_{\widetilde{M}}(2i)) = 0$
\mbox{\cite[Lecture~14]{Mu1}}. Knowing~this, the claim immediately~follows.
\item
Since~the appearance of the results of Clebsch and Coble, many mathematicians studied the moduli spaces
$\smash{\wcalM}$
and
$\calM$,
as well as the canonical morphism
$\smash{\pr\colon \wcalM \to \calM}$
connecting~them. We~do not intend to give a complete list, as this would be a hopeless~task.

But,~in addition to the references given
above, we feel that we should mention the article~\cite{CGL} of E.~Colombo, B.~van Geemen, and E.~Looijenga, where the authors reinterpret Coble's results in terms of root~systems.
For~us, some of their geometric results on the completions of the moduli spaces
$\smash{\wcalM}$
and
$\calM$
turned out to be~helpful.
\end{iii}
\end{rems}

\section{Twisting Coble's gamma variety}

\begin{ttt}\smallskip
Fix~a continuous homomorphism
$\rho\colon \Gal(\overline{K}/K) \to W(E_6)$
and consider
\begin{eqnarray*}
F_\rho \colon \{K\text{\rm -schemes}\} & \longrightarrow & \{\text{\rm sets}\} \, ,\\
S & \mapsto &\{\text{\rm marked cubic surfaces over~} S_{\overline{K}} \text{\rm ~such that~} \Gal(\overline{K}/K) \\[-1mm]
 & & \hspace{0.9cm} \text{\rm ~operates on the 27 lines as described by~} \rho \} / \!\sim \, ,
\end{eqnarray*}
the moduli functor, {\em twisted
by\/}~$\rho$.
\end{ttt}

\begin{theo}\smallskip
The functor\/
$F_\rho$
is representable by a\/
\mbox{$K$-scheme\/}
$\smash{\wcalM_{\!\!\rho}}$
that is a twist
of\/~$\smash{\wcalM}$.\medskip

\noindent
{\bf Proof.}
{\em
Let~$L/K$
be a finite Galois extension such that
$\Gal(\overline{K}/L) \subseteq \ker \rho$.
Then~the restriction
of~$F_\rho$
to the category
of~\mbox{$L$-schemes}
is clearly represented by the
\mbox{$L$-scheme}
$\smash{\wcalM_{\!L} := \wcalM \times_{\Spec K} \Spec L}$.

For~$g \in W(E_6)$,
let
$\smash{T_g\colon \wcalM_{\!L} \to \wcalM_{\!L}}$
be the morphism corresponding to the operation
of~$g$
on the 27~labels. This~is the base extension of a
morphism~$\smash{T^K_g\colon \wcalM \to \wcalM}$.
Further,~for
$\sigma \in \Gal(L/K)$,
write
$\smash{\sigma\colon \wcalM_{\!L} \to \wcalM_{\!L}}$
for the morphism induced
by~$\sigma^{-1}\colon L \leftarrow L$.
Then
\begin{eqnarray*}
\Gal(L/K) & \longrightarrow & \smash{\Mor_K(\wcalM_{\!L}, \wcalM_{\!L})} \, , \\
 \sigma & \mapsto & T_{\rho(\sigma)} \!\circ\! \sigma \, ,
\end{eqnarray*}
is a descent~datum.
Indeed,~for~$\sigma, \tau \in \Gal(L/K)$,
one~has
$$T_{\rho(\sigma)} \!\circ\! \sigma \!\circ\! T_{\rho(\tau)} \!\circ\! \tau = T_{\rho(\sigma)} \!\circ\! (\sigma \!\circ\! T_{\rho(\tau)} \!\circ\! \sigma^{-1}) \!\circ\! \sigma \!\circ\! \tau = T_{\rho(\sigma)} \!\circ\! T_{\rho(\tau)} \!\circ\! \sigma \!\circ\! \tau = T_{\rho(\sigma\tau)} \!\circ\! \sigma\tau \, .$$
Observe~that
$\sigma \!\circ\! T_{\rho(\tau)} \!\circ\! \sigma^{-1} = T_{\rho(\tau)}$,
as
$T_{\rho(\tau)}$
is the base extension of a
\mbox{$K$-morphism}.
Galois~descent \cite[Chapitre~V, \S4, n$^\circ$\,20, or~22,~Proposition~2.5]{Se1} yields a
\mbox{$K$-scheme}
$\smash{\wcalM_{\!\!\sigma}}$
such that
$\smash{\wcalM_{\!\!\sigma} \times_{\Spec K} \Spec L \cong \wcalM_{\!L}}$.

By~the universal property of the moduli
scheme~$\smash{\wcalM_{\!L}}$,
for every
\mbox{$K$-scheme}~$S$,
the set
$F_\rho(S)$
is in bijection with the set of all morphisms
$\smash{S_L \to \wcalM_{\!L}}$
of~\mbox{$L$-schemes}
such that, for every
$\sigma \in \Gal(L/K)$,
the diagram
$$
\diagram
S_L \rrto \dto_\sigma & & \smash{\wcalM_{\!L}} \dto^{T_{\rho(\sigma)} \circ \sigma} \\
S_L \rrto             & & \wcalM_{\!L}
\enddiagram
$$
commutes. Galois~descent for morphisms of schemes~\cite[Proposition~2.8]{J} shows that this datum is equivalent to giving a morphism
$\smash{S \to \wcalM_{\!\!\rho}}$
of~\mbox{$K$-schemes}.
}
\eop
\end{theo}

\begin{ttt}\smallskip
This~result suggests the following strategy to construct a smooth cubic surface
$C$
over~$\bbQ$
such that the Galois
group~$\Gal(\overline\bbQ/\bbQ)$
acts upon the lines
of~$C$
via a prescribed
subgroup~$G \subseteq W(E_6)$.\smallskip

\noindent
{\bf Strategy.}
i)
First,~find a Galois extension
$L/\bbQ$
such that
$\Gal(L/\bbQ) \cong G$.
This defines the
homomorphism~$\rho$.

\begin{iii}
\addtocounter{iii}{1}
\item
Then~a
\mbox{$\bbQ$-rational}
point~$\smash{P \in \wcalM_{\!\!\rho}(\bbQ)}$
is sought~for.
\item
For~the corresponding cubic surface
$\calC_P$
over~$\bbQ$,
the Galois group
$\Gal(\overline\bbQ/\bbQ)$
operates on the 27~lines exactly as~desired.
\end{iii}\smallskip

\noindent
Unfortunately,~we do not have the universal family
over~$\smash{\wcalM_{\!\!\rho}}$
at our disposal, at least not in a sufficiently explicit~form. Thus,~given a rational point
$\smash{P \in \wcalM_{\!\!\rho}(\bbQ)}$,
only the 40 irrational invariants
$\gamma_.$
will be known and the cubic surface has to be reconstructed from this~information. But,~anyway, searching for a
\mbox{$\bbQ$-rational}
point
on~$\smash{\wcalM_{\!\!\rho}}$
will be our main~task.
\end{ttt}

\begin{rems}\smallskip
\label{expldesc}
\begin{iii}
\item
There~is the embedding
$\gamma'\colon \wcalM_{\!L} \hookrightarrow \bP_L^{79}$
and both kinds of morphisms,
$\sigma$
and~$T_{\rho(\sigma)}$,
easily extend
to~$\bP_L^{79}$.
One~has
\begin{eqnarray*}
\sigma\colon (x_0 : \ldots : x_{79}) & \mapsto & (\sigma(x_0) : \ldots :\sigma( x_{79})) \qquad\qquad\qquad {\rm and} \\
T_{\rho(\sigma)}\colon (x_0 : \ldots : x_{79}) & \mapsto & (x_{\Pi(\rho(\sigma))^{-1}(0)} : \ldots : x_{\Pi(\rho(\sigma))^{-1}(79)}) \, .
\end{eqnarray*}
In~the second formula,
$\Pi\colon W(E_6) \hookrightarrow S_{80}$
is the permutation representation on the irrational
invariants~$\pm\gamma_.$.
To~explain why the inverses are to be taken, recall that
$T_{\rho(\sigma)}$
permutes the irrational invariants, i.e.\ the~coordinates. The~element
$x_i$
is moved to
position~$\Pi(\rho(\sigma))(i)$.
Our~formula describes exactly this~procedure.
\item
To~give a
\mbox{$K$-rational}
point
on~$\smash{\wcalM_{\!\!\rho}}$
is thus equivalent to giving an
\mbox{$L$-rational}
point
$(x_0 : \ldots : x_{79})$
on~$\smash{\gamma'(\wcalM_{\!L})}$
such that
$$(\sigma(x_{\Pi(\rho(\sigma))^{-1}(0)}) : \ldots : \sigma(x_{\Pi(\rho(\sigma))^{-1}(79)})) = (x_0 : \ldots : x_{79})$$
or, equivalently,
$(\sigma(x_0) : \ldots : \sigma(x_{79})) = (x_{\Pi(\rho(\sigma))(0)} : \ldots : x_{\Pi(\rho(\sigma))(79)})$
for
every
$\sigma \in \Gal(L/K)$.
\item
The~stronger condition that
$$(\sigma(x_{\Pi(\rho(\sigma))^{-1}(0)}), \ldots, \sigma(x_{\Pi(\rho(\sigma))^{-1}(79)})) = (x_0, \ldots, x_{79})$$
for all
$\sigma \in \Gal(L/K)$
defines a descent datum for vector spaces and, hence, a
\mbox{$80$-dimensional}
$K$-vector
space
in~$L^{80}$.

Further,~the linear relations between the irrational invariants
$\pm\gamma_.$
are generated by such with coefficients
in~$K$.
In~fact, rational numbers are possible as~coefficients. Hence,~they form an
\mbox{$L$-vector}
space that is invariant under both operations, that of
$\Gal(L/K)$
and that
of~$W(E_6)$.
This~shows that the linear relations are respected by the descent~datum. Galois~descent yields a
\mbox{$10$-dimensional}
\mbox{$K$-vector}
space~$V$
in the
\mbox{$10$-dimensional}
\mbox{$L$-vector}
space defined by the linear~relations.
\item
Analogous observations hold for the space of cubic~relations. They~form a \mbox{$30$-dimensional}
\mbox{$L$-vector}
space that is closed under the operations of
$\Gal(L/K)$
and~$W(E_6)$
and, therefore, respected by the descent~datum. Descent~yields a
\mbox{$30$-dimensional}
\mbox{$K$-vector}~space.

Consequently,~the Zariski closure of
$\smash{\wcalM_{\!\!\rho}} \subset \bP(V) \cong \bP_K^9$~is
the intersection of 30
\mbox{$K$-rational}
cubic~hypersurfaces.
\end{iii}
\end{rems}\medskip

\paragraph{\it General remarks on our approach to explicit Galois descent.}

\begin{ttt}\smallskip
\begin{iii}
\item
Our approach works as soon as we are given a finite Galois
extension~$L/K$,
a subscheme
$M \subseteq \bP^N_L$,
and a
\mbox{$K$-{\em linear\/}}
operation~$T$
of~$G := \Gal(L/K)$
on~$\bP^N_L$
such that
$M$
is invariant under
$T_\sigma \!\circ\! \sigma$
for
every~$\sigma \in G$.
Linearity~means that there is given a representation
$A\colon G \to \GL_{N+1}(K)$
such that
$T_\sigma$
is defined by the
matrix~$A(\sigma)$.

In~fact, every representation of a finite group is a subrepresentation of a sum of several copies of the regular~representation.
Consequently,~$M$
allows a linearly equivalent embedding into
some~$\bP^{N'}, N' \geq N$,
such that the
$T_\sigma$
extend to
$\bP^{N'}$
as automorphisms that simply permute the coordinates according to a permutation representation
$\pi\colon G \to S_{N'+1}$.
We~prefer permutations versus matrices in the description of the theory only in order to keep notation~concise.
\item
Consider~the particular case that the Galois descent is a~twist. I.e.,~a
\mbox{$K$-scheme}~$M_K$
is given such that
$M = M_K \times_{\Spec K} \Spec L$
and the goal is to construct another
\mbox{$K$-scheme}
$M'_K$
such
that~$M'_K \times_{\Spec K} \Spec L \cong M$.

Then~the descent datum
on~$M$
is of the form
$\{T_\sigma \!\circ\! \sigma\}_{\sigma \in G}$,
where the
$T_\sigma$
are in fact base extensions of
\mbox{$K$-scheme}~automorphisms
of~$M_K$.
What~is missing in order to apply~i) is exactly a linearization of the operation
$T\colon G \to \Aut(M)$.
\item
At~least in principle, such a linearization always exists as soon as
$M_K$
is quasi-projective. Indeed,~let
$\calL \in \Pic(M_K)$
be a very ample invertible~sheaf.
Then~$G$
operates
\mbox{$\calO_{M_K}$-linearly}
on the very ample invertible sheaf
$\smash{\bigotimes\limits_{g \in G} T_g^* \calL}$.
Use~its global sections for a projective~embedding.
\end{iii}
\end{ttt}

\section{An application to the inverse Galois problem for cubic surfaces}

\paragraph{\it A general algorithm.}

\begin{algoo}[{\rm Cubic surface for a given group}{}]\smallskip
\label{Main}
\leavevmode\\
Given~a subgroup
$G \subseteq W(E_6)$
and a field such that
$\Gal(L/\bbQ) \cong G$,
this algorithm computes a smooth cubic
surface~$C$
over~$\bbQ$
such that
$\Gal(\overline\bbQ/\bbQ)$
operates upon the lines
of~$C$
via the
group~$\Gal(L/\bbQ)$.

\begin{iii}
\item
\label{eins}
Fix~a system
$\Gamma \subseteq G$
of generators
of~$G$.
For~every
$g \in \Gamma$,
store the permutation
$\Pi(g) \in S_{80}$,
which describes the operation
of~$g$
on the 80 irrational
invariants~$\pm\gamma_.$.
Further~fix, once and for ever, ten of the
$\pm\gamma_.$
that are linearly~independent. Express~the other 70 explicitly as linear combinations of these basis~vectors.
\item
For~every
$g \in \Gamma$,
determine the
$10 \times 10$-matrix
describing the operation
of~$g$
on the
\mbox{$10$-dimensional}
\mbox{$L$-vector}
space
$\langle \gamma_. \rangle$.
Use~the explicit basis, fixed in~\ref{eins}.
\item
Choose~an explicit basis of the
field~$L$
as a
$\bbQ$-vector~space.
Finally,~make explicit the isomorphism
$\rho^{-1}\colon G \to \Gal(L/\bbQ) \subseteq \Hom_\bbQ(L,L)$.
I.e.,~write down a matrix for
every~$g \in \Gamma$.
\item
Now,~the condition that
$$(\sigma(x_{\Pi(\rho(\sigma))^{-1}(0)}), \ldots, \sigma(x_{\Pi(\rho(\sigma))^{-1}(79)})) = (x_0, \ldots, x_{79})$$
for all
$\sigma \in \Gal(L/\bbQ)$
is an explicit
\mbox{$\bbQ$-linear}
system of equations in
$10 [L : \bbQ]$~variables.
In~fact, we start with
$\Gamma$
instead
of~$\Gal(L/\bbQ)$
and get
$80 [L : \bbQ] \#\Gamma$~equations.
The~result is a ten dimensional
\mbox{$\bbQ$-vector}
space~$V \subset \langle \gamma_. \rangle$,
described by an explicit~basis.
\item
Convert the 30 cubic forms defining the image of
$\smash{\gamma_L \colon \wcalM_{\!L} \hookrightarrow \bP_L^{79}}$
into terms of this basis
of~$V$.
The~result are 30 explicit cubic forms with coefficients
in~$\bbQ$.
They~describe the Zariski closure of
$\smash{\wcalM_{\!\!\rho}}$
in~$\bP(V)$.
\item\label{sechs}
Search~for a
\mbox{$\bbQ$-rational}
point on this~variety.
\item
From~the coordinates of the point found, read the 40 irrational
invariants~$\gamma_.$.
Then~use formulas~(\ref{invvsinv}) in order to calculate Clebsch's invariants
$A, \ldots, E$.
Finally,~solve the equation problem as described in \ref{sym} and Algorithm~\ref{expl}.

In~the case that \ref{sym} or Algorithm~\ref{expl} fails, return to step~\ref{sechs}.
\end{iii}
\end{algoo}

\begin{rems}\smallskip
\begin{iii}
\item
An~important implementation trick was the~following. We~do not solve the linear system of equations
in~$L^{10}$
but
in~$\calO_L^{10}$,
for
$\calO_L \subset L$
the maximal~order. The~result is then a \mbox{rank-10}
\mbox{$\bbZ$-lattice}.
Via~the Minkowski embedding, this carries a scalar~product. Thus,~it may be reduced using the LLL-algorithm~\cite{LLL}. It~turned out in practice that points of very small height occur when taking the LLL-basis for a projective coordinate~system.

Applying~the LLL-algorithm to the lattice constructed from the maximal order should be considered as a first step towards a multivariate polynomial reduction and minimization algorithm for non-complete~intersections.
\item
There~are two points, where Algorithm~\ref{Main} may possibly~fail. First,~it may happen that no
{$\bbQ$-rational}
point is found
on~$\smash{\wcalM_{\!\!\rho}}$.
Then~one has to start with a different field having the same Galois~group.

Second,~\ref{sym} or Algorithm~\ref{expl} may fail, because of
$E = 0$,
$\Delta = 0$,
or~$F = 0$,
cf.~Remarks \ref{tec}.ii) and~iii). This~means that the cubic surface found either has no proper pentahedron, or is singular, or has nontrivial~automorphisms.

These~cases exclude a divisor from the compactified moduli space
$\bP(1,2,3,4,5)$.
Thus,~Algorithm~\ref{Main} works~generically. In~our experiments to construct examples for the remaining conjugacy classes, we met the situation that
$\Delta = 0$,
but not the situations that
$E = 0$
or~$F = 0$.
\item
In~order to get number fields with a prescribed Galois group, we used J.~Kl\"uners' number field data base {\tt http://galoisdb.math.upb.de}\,.
\end{iii}
\end{rems}\medskip

\paragraph{\it The 51 remaining conjugacy classes.}

\begin{remo}[{\rm Previous examples}{}]\smallskip
There~are exactly 350 conjugacy classes of subgroups
in~$W(E_6)$.
For~a generic cubic surface, the full
$W(E_6)$
acts upon the~lines. In~previous articles, we presented constructions producing examples for the index two subgroup
$D^1 W(E_6)$~\cite{EJ1},
all subgroups stabilizing a double-six~\cite{EJ2}, all subgroups stabilizing a pair of Steiner trihedra~\cite{EJ3}, and all subgroups stabilizing a line~\cite{EJ4}.

There~are 158 conjugacy classes stabilizing a double-six, 63 conjugacy classes stabilizing a pair of Steiner trihedra but no double-six, and 76 conjugacy classes stabilizing a line but neither a double-six nor a pair of Steiner~trihedra. Summing~up, the previous constructions completed 299 of the 350 conjugacy classes of~subgroups.
\end{remo}

\begin{ttt}\smallskip
For~some of the 51 conjugacy classes not yet covered, cubic surfaces are easily~constructed. In~fact,

\begin{iii}
\item
there~are the twists of the diagonal surface
$$X_0^3 + X_1^3 + X_2^3 + X_3^3 = 0 \, .$$
These~cubic surfaces may be written as
$\Tr_{A/\bbQ} al^3 = 0$
for
$A$
an \'etale algebra of
degree~$4$
over~$\bbQ$,
$a \in A$,
and
$l$
a linear form in four variables
over~$A$.
They~have 18 Eckardt points~\cite[section~9.1.4]{Do}.

This~approach yields nine of the 51 remaining conjugacy~classes. Their~numbers in the list are 245, 246, 289, 301, 303, 327, 337, 338, and~346.
\item
The surfaces of the type
$$\lambda X_0^3 = F_3(X_1, X_2, X_3)$$
generically have nine Eckardt points, the nine inflection points of the cubic curve, given
by~$F_3(X_1, X_2, X_3) = 0$.
This~approach yields another seven conjugacy~classes. Their~numbers are 172, 235, 236, 299, 317, 332, and~345.
\end{iii}
\end{ttt}

\begin{rem}\smallskip
In~these cases, the sets of Eckardt points are Galois~invariant. Hence,~these two constructions produce Galois groups that are contained in the
stabilizers of these~sets. These~are the two maximal subgroups of
index~$40$.
On~the other hand, the field of definition of the 27 lines
contains~$\zeta_3$,
essentially due to the Weil pairing on the relevant elliptic~curve. Thus,~there is no hope to construct in this way examples for all the groups contained in these two maximal~subgroups.
\end{rem}

\begin{ttt}\smallskip
Further,~there are a few obvious ways to try a computational brute force~attack.

\begin{iii}
\item
\label{leins}
We~systematically searched through the cubic surfaces such that all 20 coefficients are in the range
$\{-1,0,1\}$.
This~led to examples for 14 more conjugacy~classes. They~correspond to the numbers
144, 232, 267, 269, 272, 273, 305, 307, 309, 310, 329, 333, 334, 339 in the~list.
\item
Similarly,~but less systematically, we searched for cubic surfaces with a rational tritangent
plane but no rational~line. This~means, to choose a cubic field extension
$K/\bbQ$
with splitting field of type
$A_3$
or~$S_3$,
to fix a linear form
$l \in K[X_1, X_2, X_3]$,
and to search for surfaces of the~type
$$\N_{K/\bbQ} l + X_0 F_2 (X_0, X_1, X_2, X_3) = 0 \, .$$
As~there are only ten unknown coefficients, we could search in an a little bit wider~range. Note~that the generic case of this construction gives the remaining maximal subgroup of
index~$45$
in~$W(E_6)$.

Six~surfaces with orbit structures of types
$[3,12,12]$
and~$[3,24]$
have been found. The~corresponding {\tt gap}~numbers are
$90$,
$153$,
$260$,
$324$,
$335$,
and
$344$.
\item
In~analogy with~\ref{leins}, we searched through all pentahedral equations with small~coefficients. As~this family has only 5 parameters, we could inspect all surfaces with coefficients up
to~$500$.
Similarly,~we inspected all pentahedral equations with unit fractions as coefficients and denominator not more
than~$500$.
This~was motivated by simplifications shown in \cite[Fact~2.8]{EJ1b}.

This~approach results in examples for group \textnumero\,149 of order
$24$
and \textnumero\,326 of
order~$324$.
The pentahedral coefficients are
$[\frac1{256}, \frac1{241}, \frac1{225}, \frac1{81}, \frac1{81}]$
and
$[\frac1{84}, \frac1{64}, \frac1{52}, \frac1{49}, 1]$.
\item
Following~the same path, we systematically searched through all invariant vectors
$[A, \ldots, E]$
such that
$|A|, \ldots, |E| < 100$.
In~each case, we solved the equation problem as described in \ref{sym} and Algorithm~\ref{expl}.
This~led to examples for six more conjugacy~classes. Their~{\tt gap} numbers are
$216$,
$239$,
$302$,
$313$,
$319$,
and~$336$.
\end{iii}
\end{ttt}

\begin{remo}[{\rm concerning approach i)}{}]\smallskip
A priori, the search through the surfaces with small coefficients, as described in~\ref{leins}, requires the inspection of more than
$3 \!\cdot\! 10^9$~surfaces.
However,~using symmetry, we can do much better. For~this, one has to enumerate the
$3^{12}$
possible combinations of monomials of the form
$X_0^2X_1, \ldots, X_2X_3^2$.
Then~one may split this set into orbits under the operation of
$(\bbZ/2\bbZ)^4 \rtimes S_4$,
where
$S_4$
permutes the four indeterminates and
$(\bbZ/2\bbZ)^4$
changes their~signs.

This~leads to 1764~representatives. Each~representative can be extended to a
cubic surface in
$3^8$~ways by choosing coefficients for the monomials
$X_0^3$,
$X_1^3$,
$X_2^3$,
$X_3^3$,
$X_0X_1\!X_2$,
$X_0X_1\!X_3$,
$X_0X_2X_3$,
and~$X_1X_2X_3$.
Thus,~approximately
$1.1 \!\cdot\! 10^7$
surfaces had to be~inspected.
\end{remo}

\begin{remo}[{\rm concerning approaches iii) and~iv)}{}]\smallskip
Before~trying approa\-ches iii) and~iv), exactly 15 conjugacy classes were left~open. It~turned out that all these were either even, i.e.~contained in the \mbox{index-$2$}
subgroup
$D^1 W(E_6) \subset W(E_6)$,
or had a factor commutator group that was cyclic of order
$4$
or~$8$.
This~implies strong restrictions on the
discriminant~$\Delta$
of the cubic surfaces sought~for.

To~understand this, recall the following property, which partly characterizes the
discriminant~$\Delta$.
If~the 27~lines on a a smooth cubic surface
$C$
over~$\bbQ$
are acted upon by an odd Galois
group~$G \subseteq W(E_6)$
then the quadratic number field corresponding to the subgroup
$G \cap D^1 W(E_6) \subset G$
is exactly
$\smash{\bbQ\big(\!\sqrt{(-3)\Delta}\big)}$~\cite[Theorem~2.12]{EJ1}.
Correspondingly,~if
$G \subseteq W(E_6)$
is even then
$(-3)\Delta$
must be a perfect~square.

In~the odd case, the factor commutator group
$G / D^1G$
of~$G$
surjects onto
$G/G \cap D^1 W(E_6) \cong \bbZ/2\bbZ$.
Hence,~$G / D^1G$
corresponds to a
subfield~$L$
of the field of definition of the 27~lines
containing~$\smash{\bbQ\big(\!\sqrt{(-3)\Delta}\big)}$.

In~other words, there is an embedding
$\smash{\bbQ\big(\!\sqrt{(-3)\Delta}\big) \subset L}$
into a field
$L$
that is Galois and cyclic of degree of degree
$4$
(or
even~$8$)
over~$\bbQ$.
\end{remo}

\begin{lem}\smallskip
\begin{iii}
\item
If~a quadratic number field\/
$\bbQ(\sqrt{D})$
allows an embedding into a
field\/~$L$
that is Galois and cyclic of
degree\/~$4$
over\/~$\bbQ$
then\/
$D > 0$
and all prime factors\/
$p \equiv 3 \pmodulo 4$
in\/~$D$
have an even~exponent.
\item
If\/~$\bbQ(\sqrt{D})$
even allows an embedding into a field Galois and cyclic of
degree\/~$8$
then the same is true for all primes\/
$p \equiv 5 \pmodulo 8$.
\end{iii}\smallskip

\noindent
{\bf Proof.}
{\em
i) is shown in~\cite[Theorem~1.2.4]{Se2}. For~ii), the proof is~analogous. Both~results are direct applications of class field~theory.
}
\eop
\end{lem}\smallskip

\noindent
We~used this restriction in approaches iii) and~iv) as a highly efficient~pretest. It~immediately ruled out most of the~candidates.

\begin{ttt}\smallskip
To~summarize, using relatively naive methods, we found examples for 44 of the 51 remaining conjugacy~classes. Thus,~only for the last seven, we had to use the main~algorithm. In~the list, they correspond to the numbers 73, 155, 169, 177, 179, 266, 286.
\end{ttt}\medskip

\paragraph{\it Remarks concerning the running times.}

\begin{ttt}\smallskip
We~implemented the main algorithm and the elementary algorithms described in the appendix in {\tt magma}, version~2.18. We~worked on one core of an
Intel${}^{\text{(R)}}$Core${}^{\text{(TM)}}$2
Duo E8300~processor.

\begin{iii}\smallskip
\item
To~compute the numerical invariants of the gamma
variety~$\smash{\widetilde{M}}$,
given in Lemma \ref{CastMum}, the running times were less than 0.1~seconds.
\item
To~determine the coefficients in Proposition~\ref{Formeln}.ii), the running time was around 10~seconds per knot.

There~are certainly faster methods to compute the Clebsch's invariants for a given cubic~surface. We~preferred the approach described as it does not depend on deep theory and leads to compact~code. In~fact, we do much more than just calculating Clebsch's invariants, as we completely determine the~pentahedron.
\item
Our~code implementing the main algorithm for the subgroup \textnumero\,73, which is cyclic of order nine, is available on both author's web~pages as a file named {\tt c9\_example.m}. It~runs within a few seconds on the {\tt magma} online~calculator.

As~one might expect, it takes longer to run examples that involve larger number~fields. Further,~for the point search, a completely naive
\mbox{$O(N^{10})$-algorithm}
is~used. Thus,~the existence of a point of very small height is absolutely necessary for our implementation to~succeed.
\end{iii}
\end{ttt}

\begin{appendix}
\section{Some elementary algorithms}

\paragraph{\it Computing an equation from six blow-up points.}

\begin{ttt}\smallskip
\label{meth1}
Given six points
$p_1, \ldots, p_6 \in \bP^2(K)$
in general position, it is pure linear algebra to compute a sequence of 20~coefficients for the corresponding cubic~surface. First,~one has to determine a base of the kernel of a
\mbox{$6 \times 10$-matrix}
in order to find four linearly independent cubic
forms~$F_1, \ldots, F_4$
vanishing
in~$p_1, \ldots, p_6$.
To~find the cubic relation between
$F_1, \ldots, F_4$
means to solve a highly overdetermined homogeneous linear system of 220~equations in 20~variables.
\end{ttt}

\begin{rem}\smallskip
Actually,~there is a second algorithm, which is simpler but certainly less~standard. Starting~with the six points
$p_1, \ldots, p_6 \in \bP^2(K)$,
one may use formula~(85) of A.\,B.~Coble~\cite{Co1} to find hexahedral coefficients
$a_0, \ldots, a_5 \in K$
for the corresponding cubic~surface. From~this, an explicit equation is immediately~obtained.
\end{rem}\medskip

\paragraph{\it Computing the pentahedron and Clebsch's invariants from an equation.}

\begin{ttt}\smallskip
For~a cubic surface in pentahedral form,
$$C(X_0, X_1, X_2, X_3) := a_0X_0^3 + a_1X_1^3 + a_2X_2^3 + a_3X_3^3 - a_4(X_0 + X_1 + X_2 + X_3)^3 = 0$$
such that
$a_0, \ldots a_4 \in K \!\setminus\! \{0\}$,
its Hessian
$\smash{\det \frac{\partial^2 C}{\partial \!X_i \partial \!X_j} (X_0, X_1, X_2, X_3) = 0}$
has exactly ten singular~points. These~are simply the intersection points of three of the five planes defined
by~$X_0 = 0$, \ldots, $X_3 = 0$
and
$X_4 := -(X_0 + X_1 + X_2 + X_3) = 0$.
Thus,~each plane contains six of the ten singular~points.

Hence,~given a cubic surface in the form of a sequence of 20~coefficients, one has to compute its Hessian~first. If~the singular points have a configuration different from what was described then there is no~pentahedron. Otherwise,~one has to determine the five planes through six singular points and to normalize the corresponding linear forms
$l_0, \ldots, l_4$
such that their sum is~zero. To~find the five coefficients
$a_0, \ldots, a_4$
means to solve an overdetermined homogeneous linear system of 20~equations in five~variables.

There~is, however, one serious practical~difficulty. The~pentahedron is typically defined only over an
\mbox{$S_5$-ex}\-ten\-sion
of the base
field~$K$.
For~this situation, we have the following~algorithm.
\end{ttt}

\begin{algoo}[{\rm Pentahedron from cubic surface}{}]\smallskip
\label{meth2}
Let~a cubic surface
$C$
be given as a sequence of 20~coefficients. Suppose~that there is a proper pentahedron and that its field of definition is an
$S_5$-
or
\mbox{$A_5$-extension}
of the base
field~$K$.
Then~this algorithm computes the pentahedral~form.

\begin{iii}
\item
Determine~a Gr\"obner basis for the ideal
$\smash{\calI_{H_\sing} \subset K[X_0, \ldots, X_3]}$
of the singular locus of the
Hessian~$H$
of~$C$.
In~particular, this yields a univariate
\mbox{degree-$10$}
polynomial~$\overline{F}$
defining the
$S_5$-
or
\mbox{$A_5$-extension}.
\item
Uncover~a
degree-$5$
polynomial~$F$
with the same splitting~field.
When~$K = \bbQ$,
this may be done as~follows.
Run~a variant of Stauduhar's algorithm~\cite{St}. This~yields
\mbox{$p$-adic}
approximations of the ten zeroes
of~$\overline{F}$
together with an explicit description of the operation of
$S_5$
or~$A_5$.
Then~calculate
\mbox{$p$-adically}
a relative resolvent polynomial~\cite[Theorem~4]{St}, corresponding to the inclusion
$S_4 \subset S_5$
or
$A_4 \subset A_5$,
respectively. From~this, the polynomial
$F \in \bbQ[T]$
is obtained by rational~recovery.

Put~$L$
to be the extension field defined
by~$F$.
Clearly,~$[L : K] = 5$.
\item\label{step3}
Factorize~$\overline{F}$
over~$L$.
Two~irreducible factors,
$\overline{F}_1$
of
degree~$4$
and
$\overline{F}_2$
of
degree~$6$,
are~found.
\item
Determine,~in a second Gr\"obner base calculation, an element of minimal degree in the ideal
$(\calI_{H_\sing}, \overline{F}_2) \subset L[X_0, \ldots, X_3]$.
The~result is a linear
polynomial~$l$.
Its~conjugates define the five individual planes that form the~pentahedron.
\item\label{step5}
Scale~$l$
by a suitable non-zero factor
from~$L$
such
that~$\Tr_{L/K} l = 0$.
This~amounts to solving
over~$K$
a homogeneous system of four linear equations in five~variables. Then~calculate
$a \in L$
such that the equation of the surface is exactly
$\Tr_{L/K} al^3 = 0$.

Return~$a$.
Its~five conjugates are the pentahedral coefficients
of~$C$.
One~might want to return
$l$
as a second~value.
\end{iii}
\end{algoo}

\begin{rems}\smallskip
\begin{iii}
\item
Observe~that it is not necessary to perform any computations in the Galois hull
$\smash{\widetilde{L}}$
of~$L$.
\item
Let~us explain the idea behind Algorithm~\ref{meth2}. The~Galois group
$\smash{\Gal(\widetilde{L}/K) \cong S_5}$
or~$A_5$
permutes the five planes of the~pentahedron. The~ten singular points of the Hessian are in bijection with sets of three planes and permuted~accordingly.
Further,~$\smash{\Gal(\widetilde{L}/L)}$
is the stabilizer of one~plane. Under~this group, the six singular points that lie on that plane form an orbit and the four others form~another.

The~same is still true after projection to the
\mbox{$(X_0,X_1)$-line}.
Indeed,~the Galois operation immediately carries over to the~coordinates. Further,~no two of the ten points may coincide after projection, as this would define a nontrivial block structure for the image
of~$\smash{\Gal(\widetilde{L}/K)}$
in~$S_{10}$.
Our~assumptions ensure, however, that this subgroup is~primitive. This~explains the type of factorization described in step~\ref{step3}.

In~addition,
$\smash{(\calI_{H_\sing}, \overline{F}_2)}$
is the ideal of the six singular points lying on the
\mbox{$L$-ra}\-tio\-nal~plane.
That~is why a Gr\"obner base calculation for this ideal may discover the equation for that~plane.
%
\item
It~is not necessary to check the assumptions of this algorithm in~advance, as its output may be verified by a direct~calculation. Actually,~when there is no proper pentahedron, the algorithm should usually fail in the very first step, detecting that
$\smash{K[X_0, \ldots, X_3]/\calI_{H_\sing}}$
is not of length~ten. If~the Galois group is too small then more than two irreducible factors or even multiple factors may occur in step~\ref{step3}.
\item
It~would certainly be possible to make Algorithm~\ref{meth2} work for an arbitrary subgroup
of~$S_5$.
Somewhat~paradoxically, for small subgroups, the algorithm should be of lower complexity but harder to~describe. We~did not work out the details, since the present version turned out to be sufficient for our~purposes.
\item
To~compute the pentahedron for a cubic surface given by an explicit equation was considered as being a hopeless task before the formation of modern computer~algebra. The~reader might compare the concluding remarks of~\cite[section~6.6.2]{Ke}.
\end{iii}
\end{rems}

\begin{ttt}[{\rm Clebsch's invariants from pentahedral coefficients}{}]\smallskip
\label{meth3}
\leavevmode\\
Having~found the pentahedral coefficients, Clebsch's invariants may be directly calculated using formulas~(\ref{Clebsch}).
\end{ttt}

\begin{rem}\smallskip
The algorithms described up to this point were used in the proof of Proposition~\ref{Formeln}.ii).
\end{rem}\medskip

\paragraph{\it Computing a cubic surface from Clebsch's invariants. The equation problem.}

\begin{ttt}\smallskip
\label{sym}
The~other way round, given Clebsch's invariants
$[A, B, C, D, E]$
such that
$E \neq 0$,
one can calculate the corresponding base point in the pentahedral family as~follows.

Replace~$[A, B, C, D, E]$
by
$[A'\!, B'\!, C'\!, D'\!, E'] := [AE^3\!, BE^6\!, CE^9\!, DE^{12}\!, E^{16}]$
and set
$\sigma_5 := E^2$,~first.
Then~put
$\smash{\sigma_1 := \frac{B'}{\sigma_5^3}}$,
$\smash{\sigma_2 := \frac{D'}{\sigma_5^6}}$,
$\smash{\sigma_4 := \frac{C'}{\sigma_5^4}}$,
and,
finally,~$\smash{\sigma_3 := \frac{\sigma_4^2 - A'}{4\sigma_5}}$.
This~may be simplified~to
$$\textstyle [\sigma_1, \ldots, \sigma_5] = [B,D,\frac{C^2-AE}4,CE,E^2] \, .$$
\end{ttt}

\begin{rem}\smallskip
If~$\sigma_1, \ldots, \sigma_5 \in K$
then one would strongly expect that the corresponding cubic surface is defined
over~$K$.
We~learn, however, from formulas~(\ref{pent}) that
$\calC_{(\sigma_1, \ldots, \sigma_5)}$
is a priori defined only over the splitting
field~$L$
of the polynomial
$g(T) := T^5 - \sigma_1 T^4 \pm \ldots - \sigma_5 \in K[T]$.

But,~at least when
$g$
has no multiple zeroes,
$\calC_{(\sigma_1, \ldots, \sigma_5)}$
is equipped with a canonical descent~datum. Indeed,~let
$a_0, \ldots, a_4 \in L$
be the zeroes
of~$g$.
For~$\sigma \in \Gal(L/K)$,
denote by
$\pi(\sigma) \in S_5$
the corresponding permutation
of~$a_0, \ldots, a_4$.
I.e.,~$a_{\pi(\sigma)(i)} = \sigma(a_i)$.
Then~put
\begin{eqnarray*}
\Gal(L/K) & \longrightarrow & \Mor_K(\calC_{(\sigma_1, \ldots, \sigma_5)}, \calC_{(\sigma_1, \ldots, \sigma_5)}) \, , \\
\sigma & \mapsto & ((x_0 : \ldots : x_4) \mapsto (\sigma(x_{\pi(\sigma)^{-1}(0)}) : \ldots : \sigma(x_{\pi(\sigma)^{-1}(4)}))) \, .
\end{eqnarray*}
It~is easily checked that these morphisms indeed map
$\calC_{(\sigma_1, \ldots, \sigma_5)}$
onto itself and that they form a group~operation.
\end{rem}

\begin{algoo}[{\rm Computation of the Galois descent}{}]\smallskip
\label{expl}
\leavevmode\\
Given~a separable poly\-no\-mial
$g(T) = T^5 - \sigma_1 T^4 \pm \ldots - \sigma_5 \in K[T]$
of degree five, this algorithm computes the Galois descent
to~$K$
of the cubic
surface~$\calC_{(\sigma_1, \ldots, \sigma_5)}$.

\begin{iii}
\item
The~polynomial
$g$
defines an \'etale
\mbox{$K$-algebra}
$A := K[T]/(g)$.
Compute,~according to the definition, the traces
$t_i := \tr_{A/K} T^i$
for~$i = 0, \ldots, 4$.
\item
Determine the kernel of the
$1 \!\times\! 5$-matrix
$$\left(
\begin{array}{cccccc}
t_0 & t_1 & t_2 & t_3 & t_4  \\
\end{array}
\right) .$$
Choose~linearly independent kernel vectors
$(c^{(0)}_i, \ldots, c^{(4)}_i) \in K^5$
for
$i = 0, \ldots, 3$.
\item
Compute~the term
$$T \cdot \bigg[ \sum_{j=0}^4 (c_0^{(j)} X_0 + \ldots + c_3^{(j)} X_3)T^j \bigg]^3$$
modulo~$g(T)$.
This~is a cubic form
in~$X_0, \ldots, X_3$
with coefficients
in~$A$.
\item
Finally,~apply the trace coefficient-wise and output the resulting cubic form
in~$x_0, \ldots, x_3$
with 20~rational~coefficients.
\end{iii}
\end{algoo}

\begin{lem}\smallskip
For\/
$g = T^5 - \sigma_1 T^4 \pm \ldots - \sigma_5 \in K[T]$
a separable polynomial, Algorithm~\ref{expl} computes a cubic surface
over\/~$K$
that is geometrically isomorphic
to\/~$\calC_{(\sigma_1, \ldots, \sigma_5)}$.\smallskip

\noindent
{\bf Proof.}
{\em
The~\'etale algebra
$A = K[T]/(g)$
has five embeddings
$i_0, \ldots, i_4\colon A \hookrightarrow \overline{K}$
into the algebraic~closure.
For~$a_0, \ldots, a_4 \in \overline{K}$
the images
of~$T$,
we substituted into the equation
$$a_0 W_0 + \ldots + a_3 W_3 + a_4 (-W_0 - \ldots - W_3) = 0$$
the linear form
$$\textstyle l_0 := C_0X_0 + \ldots + C_3X_3 = i_0 \Big( \sum\limits_{j=0}^4 c_0^{(j)} T^j \Big) X_0 + \ldots + i_0 \Big( \sum\limits_{j=0}^4 c_3^{(j)} T^j \Big) X_3$$
and
$l_1, l_2, l_3$,
three of its~conjugates.

By~construction,
$C_0, \ldots, C_3$
form a basis of the
\mbox{$K$-vector}
space
$N \subset A$
consisting of the elements of trace~zero. In~particular,
$l_4 := -l_0 - \ldots - l_3$
is indeed the fourth~conjugate.

To~show the isomorphy, we only need to ensure that
$l_0, \ldots, l_3$
are linearly independent linear~forms. This~means that the
$5 \times 4$-matrix
$(C_i^{\sigma_j})_{0\leq j \leq 4, 0\leq i \leq 3}$
is of
rank~$4$.
Extending~$\{C_0, \ldots, C_3\}$
to a base
$\{C_0, \ldots, C_4\}$
of~$L$,
it suffices to verify that the
$5 \times 5$-matrix
$(C_i^{\sigma_j})_{0\leq j \leq 4, 0\leq i \leq 4}$
has full~rank. This~is, however, independent of the choice of the base and clear
for~$C_i = T^i$.
Indeed,~we then have a Vandermonde matrix of determinant
$\smash{\pm\prod\limits_{i<j} (T^{\sigma_i} - T^{\sigma_j}) = \pm\prod\limits_{i<j} (a_i - a_j) \neq 0}$.
}
\eop
\end{lem}

\begin{rems}\smallskip
\label{tec}
\begin{iii}
\item
It~is not hard to show that Algorithm~\ref{expl} computes the descent of the cubic
surface~$\calC_{(\sigma_1, \ldots, \sigma_5)}$
according to exactly the descent data described~above. We~skip the proof as it closely follows the lines of~\cite[Theorem~6.6]{EJ2}.
\item
Algorithm~\ref{expl} fails when
$g$~has
multiple~zeroes. For~the cubic
surface~$C$,
this means that some of its pentahedral coefficients~coincide. By~\cite[Example~9.1.25]{Do}, this is equivalent to
$C$~having
an Eckardt point, which, in turn, means that
$C$~has
a nontrivial automorphism~\cite[Theorem~9.5.8]{Do}.
Further,~there is the well-known section
$F \in \Gamma(\bP(1,2,3,4,5), \calO(25))$
that vanishes exactly on the locus corresponding to the cubic surfaces having an Eckardt~point. In~pentahedral coefficients,
$F$
is given by the expression
$I_{100}^2$~\cite[section~9.4.5]{Do}.

If~$F = 0$
then we actually face an ill-posed~problem. Due~to the presence of twists, the Clebsch invariants do not determine the cubic surface up to isomorphism
over~$K$,
but only up to isomorphism over the algebraic
closure~$\overline{K}$.
Thus,~the information available to us is insufficient on principle in order to perform a Galois~descent.
\item
Observe~that, when
$E \neq 0$
and~$F \neq 0$,
the
discriminant~$\Delta$
may nevertheless~vanish. Then~the corresponding cubic surface is~singular.
\end{iii}
\end{rems}

\begin{ttt}\smallskip
It~is classically called the {\em equation problem\/}~\cite[Definition 4.1.17]{Hu} to determine an equation for the cubic surface when the invariants
$A, \ldots, E$
are~known.
If~$E \neq 0$
and~$F \neq 0$
then \ref{sym} and Algorithm~\ref{expl} together provide an algorithmic solution to the equation~problem.
\end{ttt}

\begin{rem}\smallskip
If~$F \neq 0$
but
$E = 0$
then one might start with
$E = \varepsilon^8 \in K[\varepsilon]$
(or
$E = \varepsilon$)
instead and run Algorithm~\ref{expl} over the function~field. Unfortunately,~the resulting cubic surface typically has bad reduction
at~$\varepsilon = 0$.
Thus,~one cannot specialize
$\varepsilon$
to~$0$,~naively.
An~application of J.\,Koll\'ar's polynomial minimization algorithm~\cite[in particular Proposition 6.4.2]{Ko} is necessary to find a good~model. The~reduction
at~$\varepsilon = 0$
then solves the equation~problem.
\end{rem}
\end{appendix}

\end{document}